# BINOMIAL APPROXIMATIONS OF SHORTFALL RISK FOR GAME OPTIONS[1]

By Yan Dolinsky and Yuri Kifer

*Hebrew University*

We show that the shortfall risk of binomial approximations of game (Israeli) options converges to the shortfall risk in the corresponding Black–Scholes market considering Lipschitz continuous path-dependent payoffs for both discrete- and continuous-time cases. These results are new also for usual American style options. The paper continues and extends the study of Kifer [*Ann. Appl. Probab.* **16** (2006) 984–1033] where estimates for binomial approximations of prices of game options were obtained. Our arguments rely, in particular, on strong invariance principle type approximations via the Skorokhod embedding, estimates from Kifer [*Ann. Appl. Probab.* **16** (2006) 984–1033] and the existence of optimal shortfall hedging in the discrete time established by Dolinsky and Kifer [*Stochastics* **79** (2007) 169–195].

**1. Introduction.** This paper deals with game (Israeli) options introduced in [5] sold in a standard securities market consisting of a nonrandom component $b_t$ representing the value of a savings account at time $t$ with an interest rate $r$ and of a random component $S_t$ representing the stock price at time $t$. As usual, we view $S_t$, $t > 0$ as a stochastic process on a probability space $(\Omega, \mathcal{F}, P)$ and we assume that it generates a right-continuous filtration $\{\mathcal{F}_t\}$. The setup includes also two continuous stochastic payoff processes $X_t \geq Y_t \geq 0$ adapted to the above filtration. Recall that game contingent claim (GCC) or game option is defined as a contract between the seller and the buyer of the option such that both have the right to exercise it at any time up to a maturity date (horizon) $T$. If the buyer exercises the contract at time $t$, then he receives the payment $Y_t$, but if the seller exercises

Received March 2007; revised December 2007.
[1]Supported in part by ISF Grant 130/06.

(cancels) the contract before the buyer, then the latter receives $X_t$. The difference $\Delta_t = X_t - Y_t$ is the penalty which the seller pays to the buyer for the contract cancellation. In short, if the seller will exercise at a stopping time $\sigma \le T$ and the buyer at a stopping time $\tau \le T$, then the former pays to the latter the amount $H(\sigma, \tau)$ where

$$H(s,t) = X_s \mathbb{I}_{s<t} + Y_t \mathbb{I}_{t \le s}$$

and we set $\mathbb{I}_A = 1$ if an event $A$ occurs and $\mathbb{I}_A = 0$ if not.

A hedge (for the seller) against a GCC is defined here as a pair $(\pi, \sigma)$ which consists of a self-financing strategy $\pi$ (i.e., a trading strategy with no consumption and no infusion of capital) and a stopping time $\sigma$ which is the cancellation time for the seller. A hedge is called perfect if no matter what exercise time the buyer chooses, the seller can cover his liability to the buyer (with probability 1). The option price $V^*$ is defined as the minimal initial capital which is required for a perfect hedge, that is, for any $x > V^*$ there is a perfect hedge with an initial capital $x$. Recall (see [6]) that pricing a GCC in a complete market leads to the value of a zero sum optimal stopping (Dynkin's) game with discounted payoffs $\tilde{X}_t = b_0 \frac{X_t}{b_t}$, $\tilde{Y}_t = b_0 \frac{Y_t}{b_t}$ considered under the unique martingale measure $\tilde{P} \sim P$. The stochastic process of values $V_t^\pi$ of the portfolio $\pi$ at time $t$ is called the wealth process. In this paper we allow only hedges $(\pi, \sigma)$ with self-financing strategies $\pi$ having nonnegative wealth process, calling such $\pi$ *admissible*. This corresponds to the situation when the portfolio is handled without borrowing of the capital. In real market conditions an investor (seller) may not be willing for various reasons to tie in a hedging portfolio the full initial capital required for a perfect hedge. In this case the seller is ready to accept a risk that his portfolio value at an exercise time may be less than his obligation to pay and he will need additional funds to fulfill the contract. Thus a portfolio shortfall comes into the picture and by this reason we distinguish here between hedges and perfect hedges.

In this paper we deal with a certain type of risk called the shortfall risk (cf., e.g., [1, 2, 4, 9]) which was defined for game options in [2] by

$$R(\pi, \sigma) = \sup_{\tau} E\left[\left(Q(\sigma, \tau) - b_0 \frac{V^\pi_{\sigma\wedge\tau}}{b_{\sigma\wedge\tau}}\right)^+\right]$$

where the supremum is taken over all stopping times not exceeding a horizon $T$, $Q(s,t) = \tilde{X}_s \mathbb{I}_{s<t} + \tilde{Y}_t \mathbb{I}_{t \le s}$ is the discounted payoff, and $E$ denotes the expectation with respect to the objective probability $P$. The shortfall risk for an initial capital $x$ is defined as

$$R(x) = \inf_{(\pi,\sigma)} R(\pi, \sigma)$$

where the infimum is taken over all hedges with an initial capital $x$. An investor (seller) whose initial capital $x$ is less than the option price $V^*$ still wants to compute the minimal possible shortfall risk and to find a hedge with the initial capital $x$ which minimizes or “almost” minimizes the shortfall risk. For discrete-time models we showed in [2] how to do this but for the continuous-time Black–Scholes (BS) market the problem becomes quite complicated. The Cox, Ross and Rubinstein (CRR) binomial model (see, e.g., [12]) is an efficient tool to approximate derivative securities in a BS market. In [6] it was shown under quite general assumptions for path-dependent payoff functions that the option price (for a game option) in a BS model can be approximated by a sequence of option prices in appropriate CRR $n$-step models with errors bounded by $Cn^{-1/4}(\ln n)^{3/4}$ where $C$ is a constant which can be estimated explicitly. The main goal of this paper is to show that for path-dependent payoffs satisfying the conditions of [6] and for an initial capital $x$ the shortfall risk in a BS market $R(x)$ is a limit of the shortfall risks $R_n(x)$ for the same initial capital in an appropriate sequence of CRR markets. For game options we are able to provide only a one-sided error estimate $R(x) - R_n(x) \leq Cn^{-1/4}(\ln n)^{3/4}$ where $C > 0$ is a constant, but for American ones we derive in Section 6 full error estimates. These results rely on estimates of [6] and hedge constructions for shortfall risks in the discrete time from [2] but require also substantial additional arguments to ensure convergence under constraints.

Some discrete-time approximation results without error estimates for European options with payoffs depending only on the current stock price were obtained in [3] where the authors proved a weak convergence of shortfall risk minimizing portfolios in CRR markets to the one in the BS market. For American and Israeli options the problem was not studied before. For European options in continuous-time models (see [1, 4]) it is known that under a constraint on the initial capital there exists a portfolio which minimizes the shortfall risk. Furthermore, by using the Neyman–Pearson lemma and convex duality methods, this portfolio can be found explicitly. In [9] the author proved without an explicit construction that for American options in the continuous-time BS model there exists a portfolio which minimizes the risk. The proof was based on the fact that the shortfall risk in this case is a convex functional of the wealth process while for game options the shortfall risk fails to be a convex functional of the wealth process, and so the convex analysis methods become unavailable in this case. For game options the question whether there exists a hedge which minimizes the shortfall risk in the continuous-time BS model remains open.

In [2] we proved that for a game option in the multinomial model there exists a hedge which minimizes the shortfall risk under constraint on the initial capital, and the above hedge can be computed via a dynamical programming procedure. We will use these hedges in the CRR markets in order to construct hedges in the continuous BS market which “almost” minimize the shortfall risk. Although the BS market is continuous, in practice an investor can buy stock and bond units only on a finite set of times (may be random), and so construction of the above hedges can be useful for practical applications, since (as we will see) in order to manage the corresponding portfolios it is sufficient to buy stocks and bonds only on a finite set of random times. There was no construction of such portfolio strategies before even for European options. Our main tool is the Skorokhod type embedding of sums of i.i.d. random variables into a Brownian motion with a constant drift. This tool was employed in [6] in order to obtain error estimates for approximations of option prices. We will use this embedding in order to turn optimal hedges of CRR markets into hedges in the BS market which are almost optimal. If we could show that the sequence of the above hedges converges to a hedge in some reasonable sense, then the latter hedge would minimize the shortfall risk for the BS market, but meanwhile this problem remains open.

Main results of this paper are formulated in the next section where we discuss also the Skorokhod type embedding. In Section 3 we introduce recursive formulas which enable us to compare various risks. In Section 4 we derive auxiliary estimates for risks. In Section 5 we complete the proof of main results of the paper.

**2. Preliminaries and main results.** First, we recall the setup from [6]. Denote by $M[0,t]$ the space of Borel-measurable functions on $[0,t]$ with the uniform metric $d_{0t}(\upsilon,\tilde{\upsilon}) = \sup_{0\le s\le t} |\upsilon_s - \tilde{\upsilon}_s|$. For each $t > 0$ let $F_t$ and $\Delta_t$ be nonnegative functions on $M[0,t]$ such that for some constant $L \ge 1$ and for any $t \ge s \ge 0$ and $\upsilon,\tilde{\upsilon} \in M[0,t]$,

$$|F_s(\upsilon) - F_s(\tilde{\upsilon})| + |\Delta_s(\upsilon) - \Delta_s(\tilde{\upsilon})| \le L(s+1)d_{0s}(\upsilon,\tilde{\upsilon}) \tag{2.1}$$

and

$$\begin{aligned} &|F_t(\upsilon) - F_s(\upsilon)| + |\Delta_t(\upsilon) - \Delta_s(\upsilon)| \\ &\quad \le L\left(|t-s|\left(1 + \sup_{u\in[0,t]} |\upsilon_u|\right) + \sup_{u\in[s,t]} |\upsilon_u - \upsilon_s|\right). \end{aligned} \tag{2.2}$$

By (2.1), $F_0(\upsilon) = F_0(\upsilon_0)$ and $\Delta_0(\upsilon) = \Delta_0(\upsilon_0)$ are functions of $\upsilon_0$ only and by (2.2),

$$\begin{aligned} &F_t(\upsilon) + \Delta_t(\upsilon) \\ &\quad \le F_0(\upsilon_0) + \Delta_0(\upsilon_0) + L(t+2)\left(1 + \sup_{0\le s\le t} |\upsilon_s|\right). \end{aligned} \tag{2.3}$$

Next we consider a complete probability space $(\Omega_B, \mathcal{F}^B, P^B)$ together with a standard one-dimensional continuous-in-time Brownian motion $\{B_t\}_{t=0}^{\infty}$, and the filtration $\mathcal{F}_t^B = \sigma\{B_s | s \le t\}$. A BS financial market consists of a savings account and a stock whose prices $b_t$ and $S_t^B$ at time $t$, respectively, are given by the formulas

$$b_t = b_0 e^{rt} \quad \text{and} \quad S_t^B = S_0 e^{rt + \kappa B_t^*}, \qquad b_0, S_0 > 0, \tag{2.4}$$

where

$$B_t^* = \left(\frac{\mu}{\kappa} - \frac{\kappa}{2}\right) t + B_t, \qquad t \ge 0, \tag{2.5}$$

$r$ is the interest rate, $\kappa > 0$ is called volatility and $\mu$ is another parameter. Denote by $\tilde{S}_t^B = e^{-rt} S_t^B$ the discounted stock price. We will consider a game option in the BS market with payoff processes having the form

$$Y_t = F_t(S^B) \quad \text{and} \quad X_t = G_t(S^B), \qquad t \ge 0,$$

where $G_t = F_t + \Delta_t$ with $F$ and $\Delta$ satisfying (2.1) and (2.2), and $S^B = S^B(\omega) \in M[0, \infty)$ is a random function taking the value $S_t^B = S_t^B(\omega)$ at $t \in [0, \infty)$. When considering $F_t(S^B), G_t(S^B)$ for $t < \infty$ we take the restriction of $S^B$ to the interval $[0, t]$. Denote by $T$ the horizon of our game option assuming that $T < \infty$. Recall (see, e.g., [12], Section 7.1) that a self-financing strategy $\pi$ with a (finite) horizon $T$ and an initial capital $x$ is a process $\pi = \{\pi_t\}_{t=0}^T$ of pairs $\pi_t = (\beta_t, \gamma_t)$ where $\beta_t$ and $\gamma_t$ are progressively measurable with respect to the filtration $\mathcal{F}_t^B$, $t \ge 0$, and satisfy

$$\int_0^T e^{rt} |\beta_t| \, dt < \infty \quad \text{and} \quad \int_0^T (\gamma_t S_t^B)^2 \, dt < \infty. \tag{2.6}$$

The portfolio value $V_t^\pi$ for a strategy $\pi$ at time $t \in [0, T]$ is given by

$$V_t^\pi = \beta_t b_t + \gamma_t S_t^B = x + \int_0^t \beta_u \, db_u + \int_0^t \gamma_u \, dS_u^B. \tag{2.7}$$

Denote by $\tilde{V}_t^\pi = e^{-rt} V_t^\pi$ the discounted portfolio value at time $t$. Then it is easy to see that (see, e.g., [12])

$$\tilde{V}_t^\pi = x + \int_0^t \gamma_u \, d\tilde{S}_u^B \tag{2.8}$$

and by (2.7),

$$\beta_t = \left( x + \int_0^t \gamma_u \, d\tilde{S}_u^B - \gamma_t \tilde{S}_t^B \right) \Big/ b_0. \tag{2.9}$$

Hence, the discounted portfolio value depends only on the process $\{\gamma_t\}_{t=0}^T$ and the process $\{\beta_t\}_{t=0}^T$ can be obtained by (2.9). A self-financing strategy $\pi$ is called *admissible* if $V_t^\pi \ge 0$ for all $t \in [0, T]$ and the set of such strategies with an initial capital $x$ will be denoted by $\mathcal{A}^B(x)$. Set also $\mathcal{A}^B = \bigcup_{y>0} \mathcal{A}^B(y)$. Denote by $\mathcal{T}^B$ the set of all stopping times with respect to the Brownian filtration $\mathcal{F}_t^B$, $t \ge 0$, and let $\mathcal{T}_{0T}^B$ be the set of all stopping times with values in $[0,T]$. A pair $(\pi,\sigma) \in \mathcal{A}^B(x) \times \mathcal{T}_{0T}^B$ of an *admissible* self-financing strategy $\pi$ with an initial capital $x$ and of a stopping time $\sigma$ will be called a hedge. Set

$$Q^B(s,t) = \tilde{X}_s \mathbb{I}_{s<t} + \tilde{Y}_t \mathbb{I}_{t \le s}, \tag{2.10}$$

where $\tilde{Y}_t = e^{-rt} Y_t$ and $\tilde{X}_t = e^{-rt} X_t$ are the discounted payoffs. For a hedge $(\pi,\sigma)$ the shortfall risk is given by (see [2])

$$R(\pi,\sigma) = \sup_{\tau \in \mathcal{T}_{0T}^B} E^B[(Q^B(\sigma,\tau) - \tilde{V}^\pi_{\sigma\wedge\tau})^+], \tag{2.11}$$

which is the maximal possible expectation with respect to the probability measure $P^B$ of the discounted shortfall. The shortfall risks for a portfolio $\pi \in \mathcal{A}^B$ and for an initial capital $x$ are given by

$$R(\pi) = \inf_{\sigma \in \mathcal{T}_{0T}^B} R(\pi,\sigma) \quad \text{and} \quad R(x) = \inf_{\pi \in \mathcal{A}^B(x)} R(\pi), \tag{2.12}$$

respectively. Denote by $\tilde{P}^B$ the unique martingale measure. Using standard arguments we obtain that the restriction of the $\tilde{P}^B$ to the $\sigma$-algebra $\mathcal{F}_t^B$ satisfies

$$Z_t = \frac{dP^B}{d\tilde{P}^B}\bigg| \mathcal{F}_t^B = e^{(\mu/\kappa)B_t + (1/2)(\mu/\kappa)^2 t}. \tag{2.13}$$

By [5] the game option price $V^*$ is given by

$$V^* = \inf_{\sigma \in \mathcal{T}_{0T}^B} \sup_{\tau \in \mathcal{T}_{0T}^B} \tilde{E}^B Q^B(\sigma,\tau) \tag{2.14}$$

where $\tilde{E}^B$ is the expectation with respect to $\tilde{P}^B$.

As in [6] we consider a sequence of CRR markets on a complete probability space such that for each $n = 1, 2, \ldots$ the bond prices $b_t^{(n)}$ at time $t$ are

$$b_t^{(n)} = b_0 e^{r[nt/T]T/n} = b_0(1 + r_n)^{[nt/T]}, \qquad r_n = e^{rT/n} - 1, \tag{2.15}$$

and stock prices $S_t^{(n)}$ at time $t$ are given by the formulas $S_t^{(n)} = S_0$ for $t \in [0, T/n)$ and

$$\begin{aligned} S_t^{(n)} &= S_0 \exp\left(\sum_{k=1}^{[nt/T]} \left(\frac{rT}{n} + \kappa\left(\frac{T}{n}\right)^{1/2} \xi_k\right)\right) \\ &= S_0 \prod_{k=1}^{[nt/T]} (1 + \rho_k^{(n)}) \qquad \text{if } t \ge T/n, \end{aligned} \tag{2.16}$$

where $\rho_k^{(n)} = \exp(\frac{rT}{n} + \kappa(\frac{T}{n})^{1/2}\xi_k) - 1$ and $\xi_1, \xi_2, \ldots$ are i.i.d. random variables taking values 1 and $-1$ with probabilities $p^{(n)} = (\exp((\kappa - \frac{2\mu}{\kappa})\sqrt{\frac{T}{n}}) + 1)^{-1}$ and $1-p^{(n)}=(\exp((\frac{2\mu}{\kappa}-\kappa)\sqrt{\frac{T}{n}})+1)^{-1}$, respectively. Let $P^\xi_n=\{p^{(n)},1-p^{(n)}\}^\infty$ be the corresponding product probability measure on the space of sequences $\Omega_\xi=\{-1,1\}^\infty$ and let $\tilde S_m=(1+r_n)^{-m}S_m$ be the discounted stock price. We consider $S^{(n)}=S^{(n)}(\omega)$ as a random function on $[0,T]$, so that $S^{(n)}(\omega)\in M[0,T]$ takes the value $S^{(n)}_t=S^{(n)}_t(\omega)$ at $t\in[0,T]$. Set $\mathcal{F}^\xi_k=\sigma\{\xi_1,\dots,\xi_k\}$, $\mathcal{F}^\xi=\bigcup_{k\ge 1}\mathcal{F}^\xi_k$ and denote by $\mathcal{T}^\xi_{0n}$ the set of all stopping times with respect to the filtration $\mathcal{F}^\xi_k$ with values in $\{0,1,\dots,n\}$. Let $\mathcal{A}^{\xi,n}(x)$ be the set of all *admissible* self-financing strategies with an initial capital $x$. Recall (see [12]) that a self-financing strategy $\pi$ with an initial capital $x$ and a horizon $n$ is a sequence $(\pi_1,\dots,\pi_n)$ of pairs $\pi_k=(\beta_k,\gamma_k)$ where $\beta_k,\gamma_k$ are $\mathcal{F}^\xi_{k-1}$-measurable random variables representing the number of bond and stock units, respectively, at time $k$. Thus the portfolio value $V^\pi_k$, $k=0,1,\dots,n$ is given by

$$V^\pi_0=x,\qquad V^\pi_k=\beta_k b^{(n)}_{kT/n}+\gamma_k S^{(n)}_{kT/n},\qquad 1\le k\le n. \tag{2.17}$$

Denote by $\tilde V^\pi_k=(1+r_n)^{-k}V^\pi_k$ the discounted portfolio value at time $k$. Since $\pi$ is self-financing, then

$$\beta_k b^{(n)}_{kT/n}+\gamma_k S^{(n)}_{kT/n}=\beta_{k+1}b^{(n)}_{kT/n}+\gamma_{k+1}S^{(n)}_{kT/n}, \tag{2.18}$$

and so (see [12]),

$$\tilde V^\pi_k=x+\sum_{i=0}^{k-1}\gamma_{i+1}(\tilde S^{(n)}_{(i+1)T/n}-\tilde S^{(n)}_{iT/n}). \tag{2.19}$$

Furthermore, again,

$$\beta_k=\left(x+\sum_{i=0}^{k-1}\gamma_{i+1}(\tilde S^{(n)}_{(i+1)T/n}-\tilde S^{(n)}_{iT/n})-\gamma_k\tilde S^{(n)}_{kT/n}\right)\Big/ b_0, \tag{2.20}$$

and so, as before, in order to determine a self-financing strategy it suffices to introduce a process $\{\gamma_k\}^n_{k=0}$ and to obtain the process $\{\beta_k\}^n_{k=0}$ by (2.20). We call a self-financing strategy $\pi$ *admissible* if $V^\pi_k\ge 0$ for any $k\le n$. Set also $\mathcal{A}^{\xi,n}=\bigcup_{u>0}\mathcal{A}^{\xi,n}(u)$.

Let

$$Y^{(n)}_k=F_{kT/n}(S^{(n)}),\qquad X^{(n)}_k=G_{kT/n}(S^{(n)}) \tag{2.21}$$

and

$$Q^{(n)}(s,k)=\tilde X^{(n)}_s\mathbb{I}_{s<k}+\tilde Y^{(n)}_k\mathbb{I}_{k\le s},\qquad k,s\le n, \tag{2.22}$$

where $\tilde X^{(n)}_k=(1+r_n)^{-k}X^{(n)}_k$ and $\tilde Y^{(n)}_k=(1+r_n)^{-k}Y^{(n)}_k$ are the discounted payoffs. Clearly $Y_k,X_k$ are $\mathcal{F}^\xi_k$-measurable. A hedge with an initial capital $x$ is an element in the set $\mathcal{A}^{\xi,n}(x)\times\mathcal{T}^{\xi}_{0n}$. For a hedge $(\pi,\sigma)$ the shortfall risk is given by

$$R_n(\pi,\sigma)=\max_{\tau\in\mathcal{T}^{\xi}_{0n}}E^{\xi}_n[(Q^{(n)}(\sigma,\tau)-\tilde{V}^{\pi}_{\sigma\wedge\tau})^+],\tag{2.23}$$

which is the maximal expectation with respect to the probability measure $P^{\xi}_n$ of the discounted shortfall. Observe that $\mathcal{T}^{\xi}_{0n}$ is a finite set so that we can use max in (2.23). The shortfall risk for a portfolio $\pi\in\mathcal{A}^{\xi,n}$ and for an initial capital $x$ is given by

$$R_n(\pi)=\min_{\sigma\in\mathcal{T}^{\xi}_{0n}}R_n(\pi,\sigma)\quad\text{and}\quad R_n(x)=\inf_{\pi\in\mathcal{A}^{\xi,n}(x)}R_n(\pi),\tag{2.24}$$

respectively. Let $\tilde{P}^{\xi}_n$ be a probability measure such that $\xi_1,\xi_2,\dots$ is a sequence of i.i.d. random variables taking on the values 1 and $-1$ with probabilities $\tilde{p}^{(n)}=(\exp(\kappa\sqrt{\frac{T}{n}})+1)^{-1}$ and $1-\tilde{p}^{(n)}=(\exp(-\kappa\sqrt{\frac{T}{n}})+1)^{-1}$, respectively (with respect to $\tilde{P}^{\xi}_n$). Observe that for any $n$ the process $\{\tilde{S}^{(n)}_{mT/n}\}_{m=0}^{n}$ is a martingale with respect to $\tilde{P}^{\xi}_n$, and so we conclude that $\tilde{P}^{\xi}_n$ is the unique martingale measure for the above CRR markets.

Consider an investor in the BS market whose initial capital is $x$ which is less than the option price $V^*$. In this case the investor accepts a risk since there is no perfect hedge (see [2]). The following result says that the shortfall risk $R(x)$ of a game option in the BS market can be approximated by a sequence $R_n(x)$ of shortfall risks of game options in the CRR markets defined above and it provides also a one-sided error estimate of this approximation.

THEOREM 2.1.

$$\lim_{n\to\infty}R_n(x)=R(x).\tag{2.25}$$

*Furthermore, there exists a constant $C>0$ such that for any $n>0$,*

$$R(x)\le R_n(x)+Cn^{-1/4}(\ln n)^{3/4}.\tag{2.26}$$

Relying on convexity arguments which are not available for game options, we complement for American options in Section 6 the upper bound (2.26) by a similar lower bound.

In order to compare $R(x)$ and $R_n(x)$ we will use (a trivial form of) the Skorokhod type embedding. Thus, define recursively

$$\theta^{(n)}_0=0,\qquad \theta^{(n)}_{k+1}=\inf\left\{t>\theta^{(n)}_k:|B^*_t-B^*_{\theta^{(n)}_k}|=\sqrt{\frac{T}{n}}\right\}$$

where, recall, $B^*_t = (\frac{\mu}{\kappa} - \frac{\kappa}{2})t + B_t$. Using the same arguments as in [6] we obtain that for each of the measures $P^B, \tilde{P}^B$, the sequence $\theta^{(n)}_k - \theta^{(n)}_{k-1}$, $k = 1, 2, \dots$, is a sequence of i.i.d. random variables such that $(\theta^{(n)}_{k+1} - \theta^{(n)}_k, B^*_{\theta^{(n)}_{k+1}} - B^*_{\theta^{(n)}_k})$ are independent of $\mathcal{F}^B_{\theta^{(n)}_k}$. Employing the exponential martingale $\exp((\kappa - \frac{2\mu}{\kappa})B^*_t)$ for the probability $P^B$, we obtain that $E^B \exp((\kappa - \frac{2\mu}{\kappa})B^*_{\theta^{(n)}_1}) = 1$, concluding that $B^*_{\theta^{(n)}_1} = \sqrt{\frac{T}{n}}$ or $-\sqrt{\frac{T}{n}}$ with probability $p^{(n)}$ or $1 - p^{(n)}$, respectively. Using the martingale $\tilde{S}^B_t = S_0 \exp(\kappa B^*_t)$ for the probability $\tilde{P}^B$, we obtain $\tilde{E}^B \exp(\kappa B^*_{\theta^{(n)}_1}) = 1$, and so $B^*_{\theta^{(n)}_1} = \sqrt{\frac{T}{n}}$ or $-\sqrt{\frac{T}{n}}$ with probability $\tilde{p}^{(n)}$ or $1 - \tilde{p}^{(n)}$, respectively.

A hedge $(\pi, \sigma) \in \mathcal{A}^B(x) \times \mathcal{T}^B_{0T}$ will be called $\varepsilon$-optimal if $R(\pi, \sigma) \le R(x) + \varepsilon$. For $\varepsilon = 0$ the above hedge is called an optimal hedge. Theorem 2.1 provides an approximation of the shortfall risk of a game option in the BS market by means of the shortfall risks of game options in the CRR market which becomes especially useful if we can provide also a simple description of $\varepsilon$-optimal hedges in the BS market via optimal hedges in the CRR markets. Set $\mathfrak{b}_i = B^*_{\theta^{(n)}_i} - B^*_{\theta^{(n)}_{i-1}}$, $i = 1, 2, \dots$, and following [6] introduce for each $k = 1, 2, \dots$ the finite $\sigma$-algebra $\mathcal{G}^{B,n}_k = \sigma\{\mathfrak{b}_1, \dots, \mathfrak{b}_k\}$ with $\mathcal{G}^{B,n}_0 = \{\varnothing, \Omega_B\}$ being the trivial $\sigma$-algebra. Let $\mathcal{S}^{B,n}_{0,n}$ be the set of all stopping times with respect to the filtration $\mathcal{G}^{B,n}_k$, $k = 0, 1, 2, \dots$, with values in $\{0, 1, \dots, n\}$. Observe that for any $n$ and $k \le n$ we have a natural bijection $\Pi_{n,k} : L^\infty(\mathcal{F}^\xi_k, P^\xi_n) \to L^\infty(\mathcal{G}^{B,n}_k, P^B)$ which is given by $\Pi_{n,k} Z = \tilde{Z}$ so that if $Z = f(\xi_1, \dots, \xi_k)$ for a function $f$ on $\{-1, 1\}^k$, then $\tilde{Z} = f(\sqrt{\frac{n}{T}}\mathfrak{b}_1, \dots, \sqrt{\frac{n}{T}}\mathfrak{b}_k)$. For simplicity denote $\Pi_n = \Pi_{n,n}$ and notice that if we restrict $\Pi_n$ to $\mathcal{T}^\xi_{0n}$ we obtain a bijection $\Pi_n : \mathcal{T}^\xi_{0n} \to \mathcal{S}^{B,n}_{0,n}$. In addition to the set $\mathcal{S}^{B,n}_{0,n}$ consider also the set $\mathcal{T}^{B,n}_{0,n}$ of stopping times with respect to the filtration $\{\mathcal{F}^B_{\theta^{(n)}_k}\}_{k=0}^n$ with values in $\{0, 1, \dots, n\}$. Clearly, $\mathcal{S}^{B,n}_{0,n} \subset \mathcal{T}^{B,n}_{0,n}$. Finally, we define a function $\phi_n : \mathcal{T}^\xi_{0n} \to \mathcal{T}^B_{0T}$ which maps stopping times in CRR markets to stopping times in the BS model by

$$\phi_n(\sigma) = \begin{cases} T \wedge \theta^{(n)}_{\Pi_n(\sigma)}, & \text{if } \Pi_n(\sigma) < n, \\ T, & \text{if } \Pi_n(\sigma) = n. \end{cases} \tag{2.27}$$

Let us check that $\phi_n(\sigma) \in \mathcal{T}^B_{0T}$. Indeed, for $t < T$,

$$\{\phi_n(\sigma) \le t\} = \bigcup_{k=0}^{n-1} \{\theta^{(n)}_k \le t\} \cap \{\Pi_n(\sigma) = k\} \tag{2.28}$$

and since $\{\Pi_n(\sigma)=k\}\in\mathcal{G}_k^{B,n}\subset\mathcal{F}^B_{\theta_k^{(n)}}$ the event in the right-hand side of (2.28) belongs to $\mathcal{F}_t^B$. Since $\{\phi_n(\sigma)\le T\}=\Omega_B$ we conclude that $\phi_n(\sigma)\in\mathcal{T}^B_{0T}$. For each $n$ and $x>0$ let $\mathcal{A}^{B,n}(x)$ be the set of all *admissible* self-financing strategies with an initial capital $x$ in the BS model which can be managed only on the set $\{0,\theta_1^{(n)},\dots,\theta_n^{(n)}\}$ and such that the discounted portfolio value remains constant after the moment $\theta_n^{(n)}$. Namely, if $\pi=\{(\beta_t,\gamma_t)\}_{t=0}^\infty\in\mathcal{A}^{B,n}(x)$, then $\beta_t=\beta_{\theta_k^{(n)}}$ and $\gamma_t=\gamma_{\theta_k^{(n)}}$ provided $t\in[\theta_k^{(n)},\theta_{k+1}^{(n)})$ and $k<n$ while $\gamma_t=0$ for all $t\ge\theta_n^{(n)}$ which is achieved by selling all stocks in the portfolio at the time $\theta_n^{(n)}$, buying immediately bonds for all money and doing nothing afterward. This together with (2.8) yields that for $\pi=\{(\beta_t,\gamma_t)\}_{t=0}^\infty\in\mathcal{A}^{B,n}(x)$ the corresponding discounted portfolio value is given by

$$\tilde V_t^\pi=\begin{cases}\tilde V^\pi_{\theta_k^{(n)}}+\gamma_{\theta_k^{(n)}}(\tilde S_t^B-\tilde S^B_{\theta_k^{(n)}}), & t\in[\theta_k^{(n)},\theta_{k+1}^{(n)}],\\ \tilde V^\pi_{\theta_n^{(n)}}, & t>\theta_n^{(n)}.\end{cases}\tag{2.29}$$

Next, we define a function $\psi_n:\mathcal{A}^{\xi,n}(x)\to\mathcal{A}^{B,n}(x)$ which maps *admissible* self-financing strategies in the CRR $n$-step model to the set of *admissible* self-financing strategies in the BS model which are managed on the set $\{0,\theta_1^{(n)},\dots,\theta_n^{(n)}\}$. For $\pi=\{(\beta_k,\gamma_k)\}_{k=1}^n\in\mathcal{A}^{\xi,n}(x)$ define $\psi_n(\pi)\in\mathcal{A}^{B,n}(x)$ by

$$\tilde V_t^{\psi_n(\pi)}=\begin{cases}\tilde V^{\psi_n(\pi)}_{\theta_k^{(n)}}+\Pi_{n,k}(\gamma_{k+1})(\tilde S_t^B-\tilde S^B_{\theta_k^{(n)}}), & t\in[\theta_k^{(n)},\theta_{k+1}^{(n)}],\\ \tilde V^{\psi_n(\pi)}_{\theta_n^{(n)}}, & t>\theta_n^{(n)}.\end{cases}\tag{2.30}$$

Observe that $\Pi_{n,k}(\tilde S^{(n)}_{(kT)/n})=\tilde S^B_{\theta_k^{(n)}}$ for any $k\le n$, and so we obtain from (2.19) and (2.30) that $\tilde V^{\psi_n(\pi)}_{\theta_n^{(n)}}=\Pi_n(\tilde V_n^\pi)\ge0$. Since the discounted wealth process $\tilde V_t^{\psi_n(\pi)}$ in (2.30) is a martingale and it does not change when $t\ge\theta_n^{(n)}$, we obtain that $\tilde V_t^{\psi_n(\pi)}\ge0$ for all $t$. Hence, if $\pi$ is an *admissible* portfolio, then the portfolio $\psi_n(\pi)$ is *admissible* concluding that $\psi_n(\pi)\in\mathcal{A}^{B,n}(x)$, as required. Clearly, if we restrict the portfolio $\psi_n(\pi)$ to the interval $[0,T]$ we can consider $\psi_n(\pi)$ as an element in $\mathcal{A}^B(x)$ since the discounted wealth process $\tilde V_t^\pi$ in (2.30) is a martingale and it does not change for $t\ge\theta_n^{(n)}$, whence it is nonnegative for all $t$ if it is nonnegative at $t=\theta_n^{(n)}$.

In [2] we showed that in CRR markets, for any initial capital $x$ there exists an optimal hedge which can be calculated by a dynamical programming algorithm. We will use these hedges for our sequence of CRR markets together with the correspondence maps $\phi_n$ and $\psi_n$ introduced above in order to obtain a simple representation of $\varepsilon$-optimal hedges for the BS market.

THEOREM 2.2. *For any $n$ let $(\pi_n, \sigma_n) \in \mathcal{A}^{\xi,n}(x) \times \mathcal{T}^{\xi}_{0n}$ be the optimal hedge constructed in the next section [see (3.14) and Lemma 3.3] for the corresponding CRR markets; then*

$$\lim_{n\to\infty} R(\psi_n(\pi_n), \phi_n(\sigma_n)) = R(x). \tag{2.31}$$

**3. Optimal stopping risk representation and Skorokhod embedding.** We start with an exposition of the machinery from [2] which enables us to reduce optimization of the shortfall risk to optimal stopping problems for Dynkin's games with appropriately chosen payoff processes. For any $n$ set $a_1^{(n)} = e^{\kappa\sqrt{T/n}} - 1$, $a_2^{(n)} = e^{-\kappa\sqrt{\frac{T}{n}}} - 1$ and observe that for each $m \le n$ the random variable $\frac{\tilde{S}^{(n)}_{mT/n}}{\tilde{S}^{(n)}_{(m-1)T/n}} - 1 = \exp(\kappa(\frac{T}{n})^{1/2}\xi_m) - 1$ takes on only the values $a_1^{(n)}$ and $a_2^{(n)}$. For each $y > 0$ and $n \in \mathbb{N}$ introduce the closed interval $I_n(y) = [-\frac{y}{a_1^{(n)}}, -\frac{y}{a_2^{(n)}}]$ and for $0 \le k < n$ and a given positive $\mathcal{F}^{\xi}_k$-measurable random variable $X$, define

$$\mathcal{A}^{\xi,n}_k(X) = \biggl\{Y | Y = X + \alpha\biggl(\exp\biggl(\kappa\biggl(\frac{T}{n}\biggr)^{1/2}\xi_{k+1}\biggr) - 1\biggr) \\ \text{for some } \mathcal{F}^{\xi}_k\text{-measurable } \alpha \in I_n(X)\biggr\}. \tag{3.1}$$

Notice that if $\tilde{V}^{\pi}_k = X$ and $\tilde{V}^{\pi}_{k+1} = Y$ for $\pi = \{(\beta_k, \gamma_k)\}_{k=1}^n$, then by (2.16) and (2.19), $Y = X + \alpha(\exp(\kappa(\frac{T}{n})^{1/2}\xi_{k+1}) - 1)$ where $\alpha = \gamma_{k+1}\tilde{S}^{(n)}_{(kT)/n}$ is $\mathcal{F}^{\xi}_k$-measurable. Since we allow only nonnegative portfolio values, and so $Y \ge 0$ which must be satisfied for all possible values of $(\exp(\kappa(\frac{T}{n})^{1/2}\xi_{k+1}) - 1)$, we conclude in view of independency of $\alpha$ and $\xi_{k+1}$ that $\mathcal{A}^{\xi,n}_k(X)$ is the set of all possible discounted portfolio values at the time $k+1$ provided that the discounted portfolio value at the time $k$ is $X$.

For any $n$ and $\pi \in \mathcal{A}^{\xi,n}$ define a sequence of random variables $\{W^{\pi}_k\}_{k=0}^n$ by

$$W^{\pi}_n = (\tilde{Y}^{(n)}_n - \tilde{V}^{\pi}_n)^+ \quad \text{and} \\ W^{\pi}_k = \min((\tilde{X}^{(n)}_k - \tilde{V}^{\pi}_k)^+, \max((\tilde{Y}^{(n)}_k - \tilde{V}^{\pi}_k)^+, E^{\xi}_n(W^{\pi}_{k+1}|\mathcal{F}^{\xi}_k))) \\ \text{for } k < n. \tag{3.2}$$

Applying the results for Dynkin's games from [10] for the payoff processes

$$\{(\tilde{X}^{(n)}_k - \tilde{V}^{\pi}_k)^+\}_{k=0}^n \quad \text{and} \quad \{(\tilde{Y}^{(n)}_k - \tilde{V}^{\pi}_k)^+\}_{k=0}^n$$

in place of $\{\tilde{X}^{(n)}_k\}_{k=0}^n$ and $\{\tilde{Y}^{(n)}_k\}_{k=0}^n$ as before, we obtain that

$$W^{\pi}_0 = \min_{\sigma\in\mathcal{T}^{\xi}_{0n}} \max_{\tau\in\mathcal{T}^{\xi}_{0n}} E^{\xi}_n[(Q^{(n)}(\sigma,\tau) - \tilde{V}^{\pi}_{\sigma\wedge\tau})^+] = R_n(\pi) = R_n(\pi, \sigma(\pi)), \tag{3.3}$$

where

$$\sigma(\pi) = \min\{k|(\tilde{X}^{(n)}_k - \tilde{V}^\pi_k)^+ = W^\pi_k\} \wedge n. \tag{3.4}$$

On the Brownian probability space define $S^{B,n}_t = S_0$ if $t < T/n$ and

$$S^{B,n}_t = S_0 \exp\left(\sum_{k=1}^{[nt/T]} \left(\frac{rT}{n} + \kappa \mathfrak{b}_k\right)\right) \qquad \text{if } t \in [T/n, T] \tag{3.5}$$

where, recall, $\mathfrak{b}_k = B^*_{\theta^{(n)}_k} - B^*_{\theta^{(n)}_{k-1}}$ and consider new payoff functions $Y^{B,n}_t = F_t(S^{B,n})$ and $X^{B,n}_t = G_t(S^{B,n})$. Set

$$Q^{B,n}(s,t) = \tilde{X}^{B,n}_s \mathbb{I}_{s<t} + \tilde{Y}^{B,n}_t \mathbb{I}_{t \le s}, \tag{3.6}$$

where $\tilde{Y}^{B,n}_t = e^{-rt} Y^{B,n}_t$ and $\tilde{X}^{B,n}_t = e^{-rt} X^{B,n}_t$ are the discounted payoffs. For each positive $\mathcal{F}^B_{\theta^{(n)}_k}$-measurable random variable $X$ define $\mathcal{A}^{B,n}_k(X)$ by (3.1) with $(T/n)^{1/2}\xi_{k+1}$ and $\mathcal{F}^\xi_k$ replaced by $\mathfrak{b}_{k+1}$ and $\mathcal{F}^B_{\theta^{(n)}_k}$, respectively. By (2.29) we conclude similarly to the above that $\mathcal{A}^{B,n}_k(X)$ consists of all possible discounted values at the time $\theta^{(n)}_{k+1}$ of portfolios managed only at embedding times $\{\theta^{(n)}_i\}$ with the discounted stock evolution $\tilde{S}^B_t$ provided the discounted portfolio value at the time $\theta^{(n)}_k$ is $X$.

Next, define the shortfall risk by

$$\begin{aligned} R^{B,n}(\pi,\zeta) &= \sup_{\eta \in \mathcal{T}^{B,n}_{0n}} E^B\left[\left(Q^{B,n}\left(\frac{T\zeta}{n}, \frac{T\eta}{n}\right) - \tilde{V}^\pi_{\theta^{(n)}_{\zeta\wedge\eta}}\right)^+\right], \\ R^{B,n}(\pi) &= \inf_{\zeta \in \mathcal{T}^{B,n}_{0n}} R^{B,n}(\pi,\zeta) \quad \text{and} \quad R^{B,n}(x) = \inf_{\pi \in \mathcal{A}^{B,n}(x)} R^{B,n}(\pi). \end{aligned} \tag{3.7}$$

For any $\pi \in \mathcal{A}^{B,n}$ define a sequence of random variables $\{U^\pi_k\}^n_{k=0}$,

$$\begin{aligned} U^\pi_n &= (\tilde{Y}^{B,n}_T - \tilde{V}^\pi_{\theta^{(n)}_n})^+ \quad \text{and} \\ U^\pi_k &= \min((\tilde{X}^{B,n}_{(kT)/n} - \tilde{V}^\pi_{\theta^{(n)}_k})^+, \\ &\qquad \max((\tilde{Y}^{B,n}_{(kT)/n} - \tilde{V}^\pi_{\theta^{(n)}_k})^+, E^B(U^\pi_{k+1}|\mathcal{F}^B_{\theta^{(n)}_{k+1}}))), \qquad k < n \end{aligned} \tag{3.8}$$

and a stopping time

$$\zeta(\pi) = \min\{k|(\tilde{X}^{B,n}_{(kT)/n} - \tilde{V}^\pi_{\theta^{(n)}_k})^+ = U^\pi_k\} \wedge n. \tag{3.9}$$

Again, using the results of [10] for Dynkin's games with the adapted (with respect to the filtration $\{\mathcal{F}^B_{\theta^{(n)}_k}\}^n_{k=0}$) payoff processes $\{(\tilde{X}^{B,n}_{(kT)/n} - \tilde{V}^\pi_{\theta^{(n)}_k})^+\}^n_{k=0}$

and $\{(\tilde{Y}^{B,n}_{(kT)/n} - \tilde{V}^{\pi}_{\theta^{(n)}_k})^+\}_{k=0}^{n}$, we obtain

$$U_0^\pi = R^{B,n}(\pi) = R^{B,n}(\pi, \zeta(\pi)). \tag{3.10}$$

For $k \le n$ and $x_1, \dots, x_k \in \mathbb{R}$ consider the function $\psi^{x_1,\dots,x_k} \in M[0, \frac{kT}{n}]$ given by

$$\psi^{x_1,\dots,x_k}(t) = S_0 \exp\biggl(\frac{rjT}{n} + \kappa \sum_{i=1}^{j} x_i\biggr) \qquad \text{for } t \in [jt/n, (j+1)T/n), 1 \le j \le k$$

and

$$\psi^{x_1,\dots,x_k}(0) = S_0 \qquad \text{for } t \in [0, T/n).$$

Introduce functions $f_k^n, g_k^n : \mathbb{R}^k \to \mathbb{R}$ such that for any $x_1, \dots, x_k \in \mathbb{R}$,

$$f_k^n(x_1, \dots, x_k) = (1 + r_n)^{-k} F_{(kT)/n}(\psi^{x_1,\dots,x_k}) = e^{-rkT/n} F_{(kT)/n}(\psi^{x_1,\dots,x_k}),$$

$$g_k^n(x_1, \dots, x_k) = (1 + r_n)^{-k} G_{(kT)/n}(\psi^{x_1,\dots,x_k}) = e^{-rkT/n} G_{(kT)/n}(\psi^{x_1,\dots,x_k}).$$

Observe that for the above functions,

$$\begin{aligned} \tilde{Y}^{B,n}_{(kT)/n} &= f_k^n(\mathfrak{b}_1, \dots, \mathfrak{b}_k) \quad \text{and} \\ \tilde{X}^{B,n}_{(kT)/n} &= g_k^n(\mathfrak{b}_1, \dots, \mathfrak{b}_k), \\ \tilde{Y}^{(n)}_k &= f_k^n\biggl(\sqrt{\frac{T}{n}}\xi_1, \dots, \sqrt{\frac{T}{n}}\xi_k\biggr) \quad \text{and} \\ \tilde{X}^{(n)}_k &= g_k^n\biggl(\sqrt{\frac{T}{n}}\xi_1, \dots, \sqrt{\frac{T}{n}}\xi_k\biggr). \end{aligned} \tag{3.11}$$

The following technical lemma was proved under even more general assumptions in [2], Lemma 3.3, so for its proof we refer the reader there.

LEMMA 3.1. *Let $h_1, h_2 : [0, \infty) \to \mathbb{R}$. For a fixed $n$ define a function $\psi : [0, \infty) \to \mathbb{R}$ by*

$$\psi(y) = \inf_{u \in I_n(y)} (p^{(n)} h_1(y + u a_1^{(n)}) + (1 - p^{(n)}) h_2(y + u a_2^{(n)}))$$

*[with $p^{(n)}$ defined after (2.16)]. If $h_1, h_2$ are continuous decreasing functions, then so is $\psi$.*

For each $n$ define a sequence of functions $J_k^n : [0, \infty) \times \mathbb{R}^k \to \mathbb{R}$, $k = 0, 1, \dots, n$, by the backward recursion

$J_n^n(y, u_1, u_2, \dots, u_n)$

$$= (f_n^n(u_1,\dots,u_n) - y)^+,$$

$$\begin{aligned} J_k^n&(y,u_1,\dots,u_k)\\ &= \min\Bigg((g_k^n(u_1,\dots,u_k)-y)^+,\\ &\qquad \max\Bigg((f_k^n(u_1,\dots,u_k)-y)^+,\\ &\qquad\qquad \inf_{u\in I_n(y)}\Bigg[p^{(n)}J_{k+1}^n\Bigg(y+ua_1^{(n)},u_1,\dots,u_k,\sqrt{\frac{T}{n}}\Bigg)\\ &\qquad\qquad + (1-p^{(n)})J_{k+1}^n\Bigg(y+ua_2^{(n)},u_1,\dots,u_k,-\sqrt{\frac{T}{n}}\Bigg)\Bigg]\Bigg)\Bigg)\\ &\qquad\qquad \text{for } k=n-1,n-2,\dots,0. \end{aligned} \tag{3.12}$$

Similarly in [2], these dynamical programming relations will enable us to compare shortfall risks defined in (2.24) and (3.7) since we will be able to represent both types of risks via $J_0^n$. Meanwhile we state additional properties of the functions $J_k^n$.

Lemma 3.2. *The function $J_k^n(y,u_1,\dots,u_k)$ is continuous and decreasing with respect to $y$ for any $n$, $k \le n$.*

Proof. We fix $n$ and use the backward induction in order to prove that $J_k^n(y,u_1,\dots,u_k)$ satisfies the required conditions for any $k\le n$. For $k=n$ the statement is clear. Suppose that statement holds true for $k+1$ and prove it for $k$. Fix $u_1,\dots,u_k$. Denote $h_1(y)=J_{k+1}^n(y,u_1,\dots,u_k,\sqrt{\frac{T}{n}})$, $h_2(y)=J_{k+1}^n(y,u_1,\dots,u_k,-\sqrt{\frac{T}{n}})$ and $\psi(y)=\inf_{u\in I_n(y)}[p^{(n)}h_1(y+ua_1^{(n)}) + (1-p^{(n)})\times h_2(y+ua_2^{(n)})]$. From the induction hypothesis it follows that $h_1(y)$ and $h_2(y)$ are continuous decreasing functions, and so we obtain from Lemma 3.1 that $\psi(y)$ is continuous and decreasing, as well. Observe that

$$\begin{aligned} J_k^n&(y,u_1,\dots,u_k)\\ &= \min[(g_k^n(u_1,\dots,u_k)-y)^+, \max[(f_k^n(u_1,\dots,u_k)-y)^+,\psi(y)]], \end{aligned}$$

and so $J_k^n(y,u_1,\dots,u_k)$ is a continuous and decreasing in $y$ function. □

For a given closed interval $K=[a,b]$ and a function $f: K\times\mathbb{R}^k\to\mathbb{R}$ such that $f(\cdot,v)$ is continuous for all $v\in\mathbb{R}^k$, define $\arg\min_{a\le u\le b} f(u,v) = \min\{w\in K \mid f(w,v)=\min_{\beta\in K} f(\beta,v)\}$. The last lemma enables us to define the following functions:

$$h_k^n(y,x_1,\dots,x_k)$$

$$
\begin{aligned}
&= \operatorname*{arg\,min}_{u \in I_n(y)} \biggl[ p^{(n)} J^n_{k+1}\biggl(y + u a_1^{(n)}, u_1, \dots, u_k, \sqrt{\frac{T}{n}}\biggr) \\
&\qquad + (1 - p^{(n)}) J^n_{k+1}\biggl(y + u a_2^{(n)}, u_1, \dots, u_k, -\sqrt{\frac{T}{n}}\biggr)\biggr], \\
&\qquad k < n.
\end{aligned} \tag{3.13}
$$

Let $x$ be an initial capital. For any $n$ there exists a hedge $(\pi_n, \sigma_n) \in \mathcal{A}^{\xi,n}(x) \times \mathcal{T}^{\xi}_{0n}$ such that

$$
\begin{aligned}
&\tilde{V}^{\pi_n}_0 = x \quad \text{and} \\
&\tilde{V}^{\pi_n}_{k+1} = \tilde{V}^{\pi_n}_k + h^n_k(\tilde{V}^{\pi_n}_k, e^{\kappa\sqrt{T/n}\xi_1}, \dots, e^{\kappa\sqrt{T/n}\xi_k})(e^{\kappa\sqrt{T/n}\xi_{k+1}} - 1) \\
&\qquad \text{for } k < n \text{ and } \sigma_n = \sigma(\pi_n)
\end{aligned} \tag{3.14}
$$

with $\sigma(\pi)$ defined by (3.4). From the arguments concerning $\mathcal{A}^{\xi,n}_k(X)$ at the beginning of this section it follows that $\pi_n$ is an *admissible* strategy. From the definition of $\mathcal{A}^{B,n}_k(X)$ we conclude that for each $n$ there exists a hedge $(\tilde{\pi}_n, \zeta_n) \in \mathcal{A}^{B,n}(x) \times \mathcal{T}^{B,n}_{0,n}$ such that

$$
\begin{aligned}
&\tilde{V}^{\tilde{\pi}_n}_0 = x, \\
&\tilde{V}^{\tilde{\pi}_n}_{\theta^{(n)}_{k+1}} = \tilde{V}^{\tilde{\pi}_n}_{\theta^{(n)}_k} + h^n_k(\tilde{V}^{\tilde{\pi}_n}_{\theta^{(n)}_k}, \exp(\kappa\mathfrak{b}_1), \dots, \exp(\kappa\mathfrak{b}_k))(\exp(\kappa\mathfrak{b}_k) - 1) \quad \text{and} \\
&\zeta_n = \Pi_n(\sigma_n)
\end{aligned} \tag{3.15}
$$

with $\Pi_n$ defined before (2.27). The following lemma enables us to consider all relevant processes on the Brownian probability space and to deal with stopping times with respect to the same filtration.

Lemma 3.3. *For any $n, x > 0$,*

$$
\begin{aligned}
R_n(x) &= R_n(\pi_n) = R_n(\pi_n, \sigma_n) = J^n_0(x) \\
&= R^{B,n}(\tilde{\pi}_n) = R^{B,n}(\tilde{\pi}_n, \zeta_n) = R^{B,n}(x).
\end{aligned} \tag{3.16}
$$

Proof. For fixed $n$ and $x$ we prove first that $R^{B,n}(\tilde{\pi}_n) = R^{B,n}(x) = J^{n,x}_0(x)$. Set $\Psi^\pi_k = \tilde{V}^\pi_{\theta^{(n)}_k}$ and $\Xi_k = \tilde{V}^{\tilde{\pi}_n}_{\theta^{(n)}_k}$. We claim that for each $k \le n$ and any $\pi \in \mathcal{A}^{B,n}(x)$,

$$
J^n_k(\Psi^\pi_k, \mathfrak{b}_1, \dots, \mathfrak{b}_k) \le U^\pi_k \quad \text{and} \quad J^n_k(\Xi_k, \mathfrak{b}_1, \dots, \mathfrak{b}_k) = U^{\tilde{\pi}_n}_k. \tag{3.17}
$$

Let $k \le n$ and $\pi \in \mathcal{A}^{B,n}(x)$; then by the properties of $\mathcal{A}^{B,n}_k(\Psi^\pi_k)$ there exists a $\mathcal{F}^B_{\theta^{(n)}_k}$-measurable random variable $\alpha \in I_n(\Psi^\pi_k)$ such that $\Psi^\pi_{k+1} = \Psi^\pi_k +$ $\alpha(\exp(\kappa\mathfrak{b}_{k+1})-1)$. Since $\mathfrak{b}_{k+1}$ is independent of $\mathcal{F}^B_{\theta^{(n)}_k}$ we obtain

$$\begin{aligned}A &\stackrel{\text{def}}{=} E^B(J^n_{k+1}(\Psi^\pi_{k+1},\mathfrak{b}_1,\dots,\mathfrak{b}_{k+1})|\mathcal{F}^B_{\theta^{(n)}_k})\\ &= p^{(n)}J^n_{k+1}\left(\Psi^\pi_k+\alpha a^{(n)}_1,\mathfrak{b}_1,\dots,\mathfrak{b}_k,\sqrt{\frac{T}{n}}\right)\\ &\quad+(1-p^{(n)})J^n_{k+1}\left(\Psi^\pi_k+\alpha a^{(n)}_2,\mathfrak{b}_1,\dots,\mathfrak{b}_k,-\sqrt{\frac{T}{n}}\right),\end{aligned} \tag{3.18}$$

and so

$$\begin{aligned}A \geq \inf_{\beta\in I_n(\Psi^\pi_k)}\Bigg[&p^{(n)}J^n_{k+1}\left(\Psi^\pi_k+\beta a^{(n)}_1,\mathfrak{b}_1,\dots,\mathfrak{b}_k,\sqrt{\frac{T}{n}}\right)\\ &+(1-p^{(n)})J^n_{k+1}\left(\Psi^\pi_k+\beta a^{(n)}_2,\mathfrak{b}_1,\dots,\mathfrak{b}_k,-\sqrt{\frac{T}{n}}\right)\Bigg].\end{aligned} \tag{3.19}$$

In order to prove (3.17) we will use the backward induction. For $k=n$ the relations (3.17) follow from (3.8) and (3.12). Suppose that (3.17) hold true for $k+1$ and prove them for $k$. Let $\pi\in\mathcal{A}^{B,n}(x)$; then from (3.19) and the induction hypothesis we get

$$\begin{aligned}&E^B(U^\pi_{k+1}|\mathcal{F}^B_{\theta^{(n)}_k})\\ &\quad\geq A\geq \inf_{\beta\in I_n(\Psi^\pi_k)}\Bigg[p^{(n)}J^n_{k+1}\left(\Psi_k+\beta a^{(n)}_1,\mathfrak{b}_1,\dots,\mathfrak{b}_k,\sqrt{\frac{T}{n}}\right)\\ &\qquad\qquad+(1-p^{(n)})J^n_{k+1}\left(\Psi_k+\beta a^{(n)}_2,\mathfrak{b}_1,\dots,\mathfrak{b}_k,-\sqrt{\frac{T}{n}}\right)\Bigg].\end{aligned}$$

From (3.8) and (3.12) it follows that

$$U^\pi_k\geq J^n_k(\Psi^\pi_k,\mathfrak{b}_1,\dots,\mathfrak{b}_k). \tag{3.20}$$

Set $\alpha=h^{n,x}_k(\Xi_k,\exp(\kappa\mathfrak{b}_1),\dots,\exp(\kappa\mathfrak{b}_k))$. By the induction hypothesis similarly to (3.18) we have

$$\begin{aligned}E^B(U^{\tilde\pi_n}_{k+1}|\mathcal{F}^B_{\theta^{(n)}_k}) &= E^B(J^n_{k+1}(\Xi_{k+1},\mathfrak{b}_1,\dots,\mathfrak{b}_{k+1})|\mathcal{F}^B_{\theta^{(n)}_k})\\ &= p^{(n)}J^n_{k+1}\left(\Xi_k+\alpha a^{(n)}_1,\mathfrak{b}_1,\dots,\mathfrak{b}_k,\sqrt{\frac{T}{n}}\right)\\ &\quad+(1-p^{(n)})J^n_{k+1}\left(\Xi_k+\alpha a^{(n)}_2,\mathfrak{b}_1,\dots,\mathfrak{b}_k,-\sqrt{\frac{T}{n}}\right)\stackrel{\text{def}}{=} D.\end{aligned}$$

By the definition of $\tilde{\pi}_n$ [see (3.13)–(3.15)] we derive that

$$D = \inf_{\beta\in I_n(\Xi_k)}\left[p^{(n)} J^n_{k+1}\left(\Xi_k+\beta a^{(n)}_1,\mathfrak{b}_1,\dots,\mathfrak{b}_k,\sqrt{\frac{T}{n}}\right)\right.$$
$$\left. + (1-p^{(n)})J^n_{k+1}\left(\Xi_k+\beta a^{(n)}_2,\mathfrak{b}_1,\dots,\mathfrak{b}_k,-\sqrt{\frac{T}{n}}\right)\right].$$

From (3.8) and (3.12) it follows that

$$U^{\tilde{\pi}_n}_k = J^n_k(\Xi_k,\mathfrak{b}_1,\dots,\mathfrak{b}_k). \tag{3.21}$$

Combining (3.20) and (3.21) we obtain that (3.17) holds true for any $k$, as required. From (3.17) for $k=0$ together with (3.10) it follows that for any $\pi\in\mathcal{A}^{B,n}(x)$,

$$R^{B,n}(\tilde{\pi}_n) = U^{\tilde{\pi}_n}_0 = J^n_0(x) \le U^\pi_0 = R^{B,n}(\pi).$$

Hence,

$$J^n_0(x) = R^{B,n}(\tilde{\pi}_n) = R^{B,n}(x). \tag{3.22}$$

The proof of the equality $R_n(x)=R(\pi_n)=J^n_0(x)$ is the same; just replace $\tilde{V}^\pi_k$, $\tilde{V}^{\pi_n}_k$, $\sqrt{\frac{T}{n}}\xi_i$ and $W^\pi_k$ by $\Psi^\pi_k$, $\Xi_k$, $\mathfrak{b}_i$ and $U^\pi_k$, respectively. This together with (3.3) and (3.4) gives

$$J^n_0(x) = R_n(\pi_n,\sigma_n) = R_n(\pi_n) = R_n(x). \tag{3.23}$$

Furthermore, similarly to (3.17),

$$W^{\pi_n}_k = J^n_k\left(\tilde{V}^{\pi_n}_k,\sqrt{\frac{T}{n}}\xi_1,\dots,\sqrt{\frac{T}{n}}\xi_k\right). \tag{3.24}$$

From (3.14), (3.15), (3.21) and (3.24) we obtain that for any $k\le n$,

$$\Pi_{n,k}(\tilde{V}^{\pi_n}_k) = \tilde{V}^{\tilde{\pi}_n}_k \quad\text{and}\quad \Pi_{n,k}(W^{\pi_n}_k) = U^{\tilde{\pi}_n}_k. \tag{3.25}$$

By (3.14), $\sigma_n = \min\{k|(\tilde{X}^{(n)}_k - \tilde{V}^{\pi_n}_k)^+ = W^{\pi_n}_k\}\wedge n$ and so from (3.11) and (3.15) we have $\zeta_n = \min\{k|(\tilde{X}^{B,n}_{(kT)/n} - \tilde{V}^{\tilde{\pi}_n}_{(kT)/n})^+ = U^{\tilde{\pi}_n}_k\}\wedge n$. By (3.9) and (3.10) it follows that $R^{B,n}(\tilde{\pi}_n,\zeta_n) = R^{B,n}(\tilde{\pi}_n)$ which together with (3.22) and (3.23) completes the proof of the lemma. □

## 4. Approximations and estimates.

Set

$$A = \sup_{0\le s\le T} X_s \quad\text{and}\quad A_n = \sup_{0\le s\le \theta^{(n)}_n\vee T} X_s, \qquad n\in\mathbb{N}. \tag{4.1}$$

From the exponential moment estimates (4.8) and (4.25) of [6] it follows that there exists a constant $K_1$ such that for any natural $n$ and a real $a$,

$$E^B e^{|a|\theta_n^{(n)}\vee T} \le e^{|a|K_1 T} \quad \text{and} \quad E^B \sup_{0\le t\le \theta_n^{(n)}\vee T} \exp(aB_t) \le 2e^{a^2 K_1 T}. \tag{4.2}$$

Thus, employing the Cauchy–Schwarz inequality and (2.3), we obtain that for any $p$ there exists a constant $h_p$ such that for all $n\in\mathbb{N}$,

$$E^B A_n^p \le h_p. \tag{4.3}$$

Recall (see [12]) that for any self-financing strategy the discounted portfolio process is a right-continuous supermartingale with respect to the martingale measure. Let $\mathcal{A}^{B,M}(x)\subset\mathcal{A}^B(x)$ be the subset of all *admissible* self-financing strategies such that the corresponding discounted portfolio with the initial capital $x$ is a right-continuous martingale with respect to the martingale measure $\tilde{P}^B$ and set $\mathcal{A}^{B,M}=\bigcup_{u>0}\mathcal{A}^{B,M}(u)$.

LEMMA 4.1. *There exists a constant $K_2$ such that if $\pi,\tilde\pi\in\mathcal{A}^{B,M}$ and $\tilde{E}^B|\tilde{V}_T^\pi-\tilde{V}_T^{\tilde\pi}|<\varepsilon$, then*

$$|R^B(\pi)-R^B(\tilde\pi)| \le K_2\varepsilon^{1/4}. \tag{4.4}$$

PROOF. Let $\Upsilon=\sup_{0\le t\le T}|\tilde{V}_t^\pi-\tilde{V}_t^{\tilde\pi}|$. Using the Cauchy–Schwarz inequality we obtain

$$\begin{aligned}
&|R^B(\pi)-R^B(\tilde\pi)|\\
&\quad\le \sup_{\sigma\in\mathcal{T}_{0T}^B}\sup_{\tau\in\mathcal{T}_{0T}^B} E^B|[(Q^B(\sigma,\tau)-\tilde{V}_{\sigma\wedge\tau}^\pi)^+]-[(Q^B(\sigma,\tau)-\tilde{V}_{\sigma\wedge\tau}^{\tilde\pi})]^+|\\
&\quad\le E^B(A\mathbb{I}_{\Upsilon>\sqrt\varepsilon})+\sqrt\varepsilon=\tilde{E}^B(AZ_T\mathbb{I}_{\Upsilon>\sqrt\varepsilon})+\sqrt\varepsilon\\
&\quad\le (\tilde{E}^B A^4)^{1/4}(\tilde{E}^B Z_T^4)^{1/4}(\tilde{P}^B\{\Upsilon\ge\sqrt\varepsilon\})^{1/2}+\sqrt\varepsilon.
\end{aligned}$$

From our assumptions it follows that the process $\{|\tilde{V}_t^\pi-\tilde{V}_t^{\tilde\pi}|\}_{t=0}^T$ is a right-continuous submartingale with respect to $\tilde{P}^B$, and so using the Doob–Kolmogorov inequality we see that $\tilde{P}^B\{\Upsilon\ge\sqrt\varepsilon\}\le\frac{\tilde{E}^B|\tilde{V}_T^\pi-\tilde{V}_T^{\tilde\pi}|}{\sqrt\varepsilon}\le\sqrt\varepsilon$. Thus (assuming $\varepsilon<1$) we obtain (4.4) with $K_2=(\tilde{E}^B A^4)^{1/4}\tilde{E}^B(Z_T{}^4)^{1/4}+1$, completing the proof. □

Set $\tilde{B}_t=B_t+\frac{\mu}{\kappa}t$. From Girsanov's theorem it follows that $\{\tilde{B}_t\}_{t=0}^T$ is a Brownian motion with respect to the martingale measure $\tilde{P}^B$ and the filtration $\mathcal{F}_t^B$. The following lemma is a standard result but since we could not find a direct reference its proof for the reader's convenience is given here.

Lemma 4.2. *For any nonnegative random variable* $X \in L^1(\mathcal{F}_T, \tilde{P}^B)$, $X \neq 0$ *and* $\varepsilon > 0$ *there exist* $t_1, \ldots, t_k \in [0,T]$ *and a smooth function with a compact support* $0 \le g \in C_0^\infty(\mathbb{R}^k)$ *such that*

$$\tilde{E}^B|X - g(B_{t_1}, \ldots, B_{t_k})| < \varepsilon \quad \text{and} \quad \tilde{E}^B g(B_{t_1}, \ldots, B_{t_k}) < \tilde{E}^B X. \tag{4.5}$$

Proof. Observe that without loss of generality we can assume that $\tilde{E}^B X = 1$. Fix $\varepsilon > 0$. It is well known (see Lemma 4.3.1 in [11]) that there exist $t_1, \ldots, t_k \in [0,T]$ and a smooth function with a compact support $0 \le f \in C_0^\infty(\mathbb{R}^k)$ such that $\tilde{E}^B|X - f(\tilde{B}_{t_1}, \ldots, \tilde{B}_{t_k})| < \frac{\varepsilon}{2}$. Set $h = \frac{f}{1+\varepsilon/2}$ and observe that $\tilde{E}^B h(\tilde{B}_{t_1}, \ldots, \tilde{B}_{t_k}) < \frac{\tilde{E}^B X + \varepsilon/2}{1+\varepsilon/2} = \tilde{E}^B X$. Furthermore,

$$\begin{aligned}
&\tilde{E}^B|X - h(\tilde{B}_{t_1}, \ldots, \tilde{B}_{t_k})| \\
&\quad \le \tilde{E}^B|X - f(\tilde{B}_{t_1}, \ldots, \tilde{B}_{t_k})| + \tilde{E}^B|f(\tilde{B}_{t_1}, \ldots, \tilde{B}_{t_k}) - h(\tilde{B}_{t_1}, \ldots, \tilde{B}_{t_k})| \\
&\quad \le \frac{\varepsilon}{2} + \frac{\varepsilon}{2}\tilde{E}^B h(\tilde{B}_{t_1}, \ldots, \tilde{B}_{t_k}) < \varepsilon.
\end{aligned}$$

Next define a function $g \in C_0^\infty(\mathbb{R}^k)$ by $g(x_1, \ldots, x_k) = h(x_1 + \frac{\mu}{\kappa}t_1, \ldots, x_1 + \frac{\mu}{\kappa}t_k)$ and the result follows. $\square$

For any $x$ let $\mathcal{A}^{B,C}(x) \subset \mathcal{A}^{B,M}(x)$ be the subset consisting of all $\pi \in \mathcal{A}^{B,M}(x)$ such that $\tilde{V}_T^\pi = g(B_{t_1}, \ldots, B_{t_k})$ for some smooth function $g \in C_0^\infty(\mathbb{R}^k)$ with a compact support and $t_1, \ldots, t_k \in [0,T]$.

Lemma 4.3. *For any initial capital* $x$ *and* $\varepsilon > 0$ *there exist* $y < x$ *and* $\pi \in \mathcal{A}^{B,C}(y)$ *such that*

$$R(\pi) < R(x) + \varepsilon. \tag{4.6}$$

Proof. Fix $x, \varepsilon$ and let $\tilde{\pi} \in \mathcal{A}^B(x)$ satisfy $R(\tilde{\pi}) < R(x) + \frac{\varepsilon}{2}$. Set $M_t = \tilde{V}_t^{\tilde{\pi}} \wedge D_t$ where $\{D_t\}_{t=0}^T$ is the regular continuous martingale defined by $D_t = \tilde{E}^B(A|\mathcal{F}_t^B)$ where $A$ is the same as in (4.1). Observe that under $\tilde{P}^B$, $\{M_t\}_{t=0}^T$ is a right-continuous supermartingale which belongs to the class $\mathcal{D}$ (see, e.g., [8]). Using the Doob–Meyer decomposition (see [8]) we obtain that there exists a right-continuous martingale $\{\tilde{M}_t\}_{t=0}^T$ belonging to the class $\mathcal{D}$ and a positive adapted process $\{U_t\}_{t=0}^T$ such that

$$U_0 = 0 \quad \text{and} \quad M_t = \tilde{M}_t - U_t.$$

Thus $\tilde{M}_0 = M_0 = x \wedge D_0 \le x$. Let $\delta = (\frac{\varepsilon}{2K_2})^4$ where $K_2$ is a constant from Lemma 4.1. By Lemma 4.2 we obtain that there exist $0 \le g \in C_0^\infty(\mathbb{R}^k)$ and $t_1, \ldots, t_k \in [0,T]$ such that

$$\tilde{E}^B|\tilde{M}_T - g(B_{t_1}, \ldots, B_{t_k})| < \delta \quad \text{and} \quad \tilde{E}^B g(B_{t_1}, \ldots, B_{t_k}) < \tilde{E}^B \tilde{M}_T \le x. \tag{4.7}$$

Set $y=\tilde{E}^B g(B_{t_1},\dots,B_{t_k})$. It follows from (4.7) that $y<x$. Since the BS market is complete there exists $\pi\in\mathcal{A}^{B,C}(y)$ such that $\tilde{V}^{\pi}_t=\tilde{E}^B(g(B_{t_1},\dots,B_{t_k})|\mathcal{F}_t)$. By Lemma 3.4 we obtain that

$$\begin{aligned}R(\pi)&\le K_2\delta^{1/4}+\inf_{\sigma\in\mathcal{T}^B_{0T}}\sup_{\tau\in\mathcal{T}^B_{0T}}E^B[(Q^B(\sigma,\tau)-\tilde{M}_{\sigma\wedge\tau})^+]\\ &\le K_2\delta^{1/4}+\inf_{\sigma\in\mathcal{T}^B_{0T}}\sup_{\tau\in\mathcal{T}^B_{0T}}E^B[(Q^B(\sigma,\tau)-M_{\sigma\wedge\tau})^+].\end{aligned}\tag{4.8}$$

Since $D_t\ge X_t$, then for any $\sigma,\tau\in\mathcal{T}^B_{0T}$, $(Q^B(\sigma,\tau)-M_{\sigma\wedge\tau})^+=(Q^B(\sigma,\tau)-\tilde{V}^{\tilde{\pi}}_{\sigma\wedge\tau})^+$. Hence, by (4.8),

$$R(\pi)\le K_2\delta^{1/4}+R(\tilde{\pi})<R(x)+\varepsilon,$$

completing the proof. □

Next, we prove a general result employing arguments similar to the proof of Lemma 3.2 in [6].

Lemma 4.4. *Let $n\in\mathbb{N}$ and $\tau_1,\tau_2\le\theta^{(n)}_n\vee T$ be stopping times with respect to the Brownian filtration. Then there exist constants $L_1,L_2$ such that*

$$\text{(i)}\quad E^B|e^{-r\tau_1}F_{\tau_1}(S^B)-e^{-r\tau_2}F_{\tau_2}(S^B)|\le L_1(E^B(\tau_1-\tau_2)^2)^{1/2}+L_2(E^B(\tau_1-\tau_2)^2)^{1/4}$$

*and*

$$\text{(ii)}\quad E^B|e^{-r\tau_1}G_{\tau_1}(S^B)-e^{-r\tau_2}G_{\tau_2}(S^B)|\le L_1(E^B(\tau_1-\tau_2)^2)^{1/2}+L_2(E^B(\tau_1-\tau_2)^2)^{1/4},$$

*where, recall, $F_t$ and $G_t=F_t+\Delta_t$ are functions introduced at the beginning of Section 2.*

Proof. We start with the first term. By the Cauchy–Schwarz inequality

$$\begin{aligned}&E^B|e^{-r\tau_1}F_{\tau_1}(S^B)-e^{-r\tau_2}F_{\tau_2}(S^B)|\\ &\quad\le E^B(|e^{-r\tau_1}-e^{-r\tau_2}|F_{\tau_2}(S^B))+E^B|e^{-r\tau_1}F_{\tau_1}(S^B)-e^{-r\tau_1}F_{\tau_2}(S^B)|\\ &\quad\le rE^B[|\tau_1-\tau_2|A_n]+E^B|F_{\tau_1}(S^B)-F_{\tau_2}(S^B)|\\ &\quad\le rh_2^{1/2}(E^B(\tau_1-\tau_2)^2)^{1/2}+E^B|F_{\tau_1}(S^B)-F_{\tau_2}(S^B)|\end{aligned}\tag{4.9}$$

with $h_2$ the same as in (4.3). Using (2.2) we obtain that

$$|F_{\tau_1}(S^B)-F_{\tau_2}(S^B)|\le I_1+I_2,\tag{4.10}$$

where

$$I_1 = L(\tau_1 \vee \tau_2 - \tau_1 \wedge \tau_2)\Big(1 + \sup_{0 \le t \le \theta_n^{(n)} \vee T} S_t^B\Big),$$

$$I_2 = \sup_{\tau_1 \wedge \tau_2 \le t \le \tau_1 \vee \tau_2} L|S_t^B - S_{\tau_1 \wedge \tau_2}^B|.$$

Again, using the Cauchy–Schwarz inequality and (4.2) we obtain that there exists a constant $H^{(1)}$ such that

$$E^B I_1 \le H^{(1)}(E^B(\tau_1 - \tau_2)^2)^{1/2}. \tag{4.11}$$

Observe that

$$S_t^B = S_0 + \kappa \int_0^t S_u^B \, dB_u + (r + \mu) \int_0^t S_u^B \, du. \tag{4.12}$$

Using the Doob–Kolmogorov inequality and Itô's isometry for stochastic integrals (see, e.g., [11]) we obtain

$$\begin{aligned}
E^B & \sup_{\tau_1 \wedge \tau_2 \le t \le \tau_1 \vee \tau_2} \left| \int_{\tau_1 \wedge \tau_2}^t S_u^B \, dB_u \right| \\
&\le \left( E^B \sup_{\tau_1 \wedge \tau_2 \le t \le \tau_1 \vee \tau_2} \left| \int_{\tau_1 \wedge \tau_2}^t S_u^B dB_u \right|^2 \right)^{1/2} \\
&\le 2 \left( E^B \left( \int_{\tau_1 \wedge \tau_2}^{\tau_1 \vee \tau_2} S_u^B \, dB_u \right)^2 \right)^{1/2} \\
&= 2 \left( E^B \int_{\tau_1 \wedge \tau_2}^{\tau_1 \vee \tau_2} (S_u^B)^2 \, du \right)^{1/2} \\
&\le 2 \left( E^B \Big( |\tau_1 - \tau_2| \sup_{0 \le t \le \theta_n^{(n)} \vee T} (S_t^B)^2 \Big) \right)^{1/2}.
\end{aligned}$$

This together with (4.12) and the Cauchy–Schwarz inequality yields

$$\begin{aligned}
E^B I_2 &\le 2L\kappa \left( E^B \Big( |\tau_1 - \tau_2| \sup_{0 \le t \le \theta_n^{(n)} \vee T} (S_t^B)^2 \Big) \right)^{1/2} \\
&\quad + |r + \mu| L E^B \Big( |\tau_1 - \tau_2| \sup_{0 \le t \le \theta_n^{(n)} \vee T} S_t^B \Big) \\
&\le H^{(2)}(E^B(\tau_1 - \tau_2)^2)^{1/2} + \tilde{H}^{(2)}(E^B(\tau_1 - \tau_2)^2)^{1/4}
\end{aligned}$$

for some constants $H^{(2)}, \tilde{H}^{(2)}$. Combining (4.9)–(4.11) and (4.13) we complete the proof of (i) while (ii) is derived in a same way with the same constants. □

**5. Proving the main results.** In this section we complete the proof of Theorems 2.1 and 2.2, relying on the key Lemma 3.3, on estimates and on approximation results from Section 4 and on some additional estimates similar to [6]. We start with the lower bound estimate of the BS risk where we can rely only on quite general Lemmas 4.2 and 4.3 which do not provide specific estimates and a good lower bound in Theorem 2.1 would require more precise information on optimal hedges of shortfall risk in the BS market. Concerning the upper bound estimate we observe that admissible portfolio strategies which are managed only at embedding times are also admissible portfolio strategies for the continuous BS market which will lead to the estimate (2.26).

Let $x$ be an initial capital and $\varepsilon > 0$. From Lemma 4.3 it follows that there exist $k$, $0 < t_1 < t_2 < \cdots < t_k \le T$ and $0 \le f, g \in C_0^\infty(\mathbb{R}^k)$ such that $f(x_1,\dots,x_k) = g(x_1 + \frac{\kappa}{2}t_1,\dots,x_k + \frac{\kappa}{2}t_k)$, and so $f(B^*_{t_1},\dots,B^*_{t_k}) = g(\tilde{B}_{t_1},\dots,\tilde{B}_{t_k})$ while the portfolio $\pi \in \mathcal{A}^B$ with $\tilde{V}^\pi_t = \tilde{E}(f(B^*_{t_1},\dots,B^*_{t_k})|\mathcal{F}_t)$ satisfies

$$R(\pi) < R(x) + \varepsilon \quad \text{and} \quad V^\pi_0 < x. \tag{5.1}$$

Set

$$\Psi_n = f(B^*_{\theta^{(n)}_{[nt_1/T]}},\dots,B^*_{\theta^{(n)}_{[nt_k/T]}}), \tag{5.2}$$

$u_n = \max_{0\le k\le n} |\theta^{(n)}_k - \frac{kT}{n}|$ and $w_n = \max_{0\le k\le n-1} |\theta^{(n)}_k - \theta^{(n)}_{k-1}| + |T - \theta^{(n)}_n|$. Since $w_n \le 3u_n + \frac{T}{n}$, then from (4.7) in [6] we obtain that for any $m$ there exists a constant $K^{(m)}$ such that for all $n$,

$$E^B u_n^{2m} \le K^{(m)} n^{-m} \quad \text{and} \quad E^B w_n^{2m} \le K^{(m)} n^{-m}. \tag{5.3}$$

Clearly, $(B^*_t - B^*_{\theta^{(n)}_{[nt/T]}})^2 \le 2(B_t - B_{\theta^{(n)}_{[nt/T]}})^2 + 2((\frac{\mu}{\kappa} - \frac{\kappa}{2})(t - \theta^{(n)}_{[nt/T]}))^2$ and $|t - \theta^{(n)}_{[nt/T]}| \le \frac{T}{n} + u_n$. Hence, from (5.3) and the Doob–Kolmogorov inequality it follows that there exists a constant $H^{(3)}$ such that for all $t$, $E^B|B^*_t - B^*_{\theta^{(n)}_{[nt/T]}}|^2 \le H^{(3)} n^{-1/2}$. Let $L(f) = \max_{1\le i\le k} \sup_{x\in\mathbb{R}^k} |\frac{\partial f}{\partial x_i}(x_1,\dots,x_k)|$. Then by (5.2) and the inequality $(\sum_{i=1}^k a_i)^2 \le k\sum_{i=1}^k a_i^2$ we obtain

$$\begin{aligned} E^B(\Psi_n - \tilde{V}^\pi_T)^2 &\le L(f)^2 E^B\left(\sum_{i=1}^k |B^*_{t_k} - B^*_{\theta^{(n)}_{[nt_k/T]}}|\right)^2 \\ &\le kL(f)^2 \sum_{i=1}^k E^B(B^*_{t_k} - B^*_{\theta^{(n)}_{[nt_k/T]}})^2 \le k^2 L(f)^2 H^{(3)} n^{-1/2}. \end{aligned} \tag{5.4}$$

By (4.2) and the Cauchy–Schwarz inequality,

$$\lim_{n\to\infty} \tilde{E}^B|\Psi_n - \tilde{V}^\pi_T| = \lim_{n\to\infty} (E^B|\Psi_n - \tilde{V}^\pi_T|^2)^{1/2} (E^B Z^{-2}_{\theta^{(n)}_n \vee T})^{1/2} = 0,$$

where $Z_t$ is the Radon–Nikodym derivative given by (2.13). Since $\tilde{E}^B\tilde{V}^\pi_T < x$, then for sufficiently large $n$ we can assume that $v_n = \tilde{E}(\Psi_n) < x$. Observe that the finite-dimensional distributions of the sequence $\sqrt{\frac{T}{n}}\xi_1, \dots, \sqrt{\frac{T}{n}}\xi_n$ with respect to $\tilde{P}^\xi_n$ and the finite-dimensional distributions of the sequence $B^*_{\theta^{(n)}_1}, \dots, B^*_{\theta^{(n)}_n} - B^*_{\theta^{(n)}_{n-1}}$ with respect to $\tilde{P}^B$ are the same, and so (for sufficiently large $n$),

$$v_n = \tilde{E}^\xi_n f\left(\sqrt{\frac{T}{n}}\sum_{i=1}^{[nt_1/T]}\xi_i, \dots, \sqrt{\frac{T}{n}}\sum_{i=1}^{[nt_k/T]}\xi_i\right) < x.$$

Since CRR markets are complete we can find a portfolio $\tilde{\pi}(n) \in \mathcal{A}^{\xi,n}(v_n)$ such that

$$\tilde{V}^{\tilde{\pi}}_n = f\left(\sqrt{\frac{T}{n}}\sum_{i=1}^{[nt_1/T]}\xi_i, \dots, \sqrt{\frac{T}{n}}\sum_{i=1}^{[nt_k/T]}\xi_i\right). \tag{5.5}$$

Let $\pi' = \psi_n(\tilde{\pi}) \in \mathcal{A}^{B,n}(v_n)$; then by the definition (2.30), $\tilde{V}^{\pi'}_{\theta^{(n)}_n} = \Psi_n$. Since $R_n(\cdot)$ is a nonincreasing function, then by (5.1),

$$R_n(x) - R(x) \le R_n(v_n) - R(x) \le \varepsilon + R^{B,n}(\pi') - R(\pi). \tag{5.6}$$

Given $\delta > 0$ there exists a stopping time $\sigma(\delta) \in \mathcal{T}^B_{0T}$ such that

$$R(\pi) > \sup_{\tau\in\mathcal{T}^B_{0T}} E^B[(Q^B(\sigma,\tau) - \tilde{V}^\pi_{\sigma\wedge\tau})^+] - \delta. \tag{5.7}$$

Define a stopping time $\zeta = \zeta(n,\sigma) \in \mathcal{T}^{B,n}_{0,n}$ by

$$\zeta = \begin{cases} n \wedge \min\{i | \theta^{(n)}_i \ge \sigma\}, & \text{if } \sigma < T, \\ n, & \text{if } \sigma = T. \end{cases}$$

Next, check that $\zeta \in \mathcal{T}^{B,n}_{0,n}$. Since the Brownian filtration is right-continuous, then for any $i < n$, $\{\zeta \le i\} = \{\sigma \le \theta^{(n)}_i\} \cap \{\sigma < T\} \in \mathcal{F}^B_{\theta^{(n)}_i}$ and $\{\zeta \le n\} = \Omega_B$, thus $\zeta \in \mathcal{T}^{B,n}_{0,n}$. Clearly, there exists a stopping time $\eta = \eta(n,\zeta)$ such that

$$\begin{aligned} &E^B[(Q^{B,n}(\theta^{(n)}_\zeta, \theta^{(n)}_\eta) - \tilde{V}^{\pi'}_{\theta^{(n)}_{\zeta\wedge\eta}})^+] \\ &\quad > \sup_{\tilde{\eta}\in\mathcal{T}^{B,n}_{0,n}} E^B[(Q^{B,n}(\theta^{(n)}_\zeta, \theta^{(n)}_{\tilde{\eta}}) - \tilde{V}^{\pi'}_{\theta^{(n)}_{\zeta\wedge\tilde{\eta}}})^+] - \delta \ge R^{B,n}(\pi') - \delta. \end{aligned} \tag{5.8}$$

Similarly to Lemmas 3.2 and 3.3 in [6] it follows that there exists a constant $C_1$ such that for any $n$,

$$\sup_{\zeta\in\mathcal{T}^{B,n}_{0,n}} \sup_{\eta\in\mathcal{T}^{B,n}_{0,n}} E^B\left[\left|Q^B(\theta^{(n)}_\zeta, \theta^{(n)}_\eta) - Q^{B,n}\left(\frac{\zeta T}{n}, \frac{\eta T}{n}\right)\right|\right]$$

$$\le C_1 n^{-1/4}(\ln n)^{3/4}. \tag{5.9}$$

Observe that if $\sigma \ge \theta_\eta^{(n)} \wedge T$, then $\zeta \ge \eta$, and so from (5.7)–(5.9) we obtain

$$\begin{aligned} &R^{B,n}(\pi') - R(\pi) \\ &\quad < C_1 n^{-1/4}(\ln n)^{3/4} + 2\delta + E^B|\tilde{V}^{\pi'}_{\theta^{(n)}_{\zeta\wedge\eta}} - \tilde{V}^{\pi}_{\sigma\wedge\theta^{(n)}_\eta}| \\ &\quad\quad + E^B[(Q^B(\theta^{(n)}_\zeta, \theta^{(n)}_\eta) - Q^B(\sigma, \theta^{(n)}_\eta \wedge T))^+] \\ &\quad \le J_1 + J_2 + I_1 + I_2 + 2\delta + C_1 n^{-1/4}(\ln n)^{3/4} \end{aligned} \tag{5.10}$$

where

$$I_1 = E^B|\tilde{V}^{\pi'}_{\theta^{(n)}_{\zeta\wedge\eta}} - \tilde{V}^{\pi}_{\theta^{(n)}_{\zeta\wedge\eta}\wedge T}|, \qquad I_2 = E^B|\tilde{V}^{\pi}_{\theta^{(n)}_{\zeta\wedge\eta}\wedge T} - \tilde{V}^{\pi}_{\theta^{(n)}_\eta\wedge\sigma}|$$

and since $|\theta^{(n)}_{\zeta\wedge\eta} - \theta^{(n)}_\eta \wedge \sigma| \le w_n$, then by (5.3) and Lemma 4.4,

$$\begin{aligned} J_1 &= E^B|e^{-r\theta^{(n)}_{\zeta\wedge\eta}} G_{\theta^{(n)}_{\zeta\wedge\eta}}(S^B) - e^{-r\sigma\wedge\theta^{(n)}_\eta} G_{\theta^{(n)}_\eta\wedge\sigma}(S^B)| \le H^{(4)} n^{-1/4}, \\ J_2 &= E^B|e^{-r\theta^{(n)}_{\zeta\wedge\eta}} F_{\theta^{(n)}_{\zeta\wedge\eta}}(S^B) - e^{-r\sigma\wedge\theta^{(n)}_\eta} F_{\theta^{(n)}_\eta\wedge\sigma}(S^B)| \le H^{(4)} n^{-1/4} \end{aligned} \tag{5.11}$$

for some constant $H^{(4)}$. Clearly,

$$\begin{aligned} \tilde{V}^{\pi'}_{\theta^{(n)}_{\zeta\wedge\eta}} - \tilde{V}^{\pi}_{\theta^{(n)}_{\zeta\wedge\eta}\wedge T} &= \tilde{E}^B(\Psi_n - \tilde{V}^\pi_T | \mathcal{F}_{\theta^{(n)}_{\zeta\wedge\eta}}) \\ &= E^B\Big(\frac{Z_{\theta^{(n)}_{\zeta\wedge\eta}}}{Z_{T\vee\theta^{(n)}_n}}(\Psi_n - \tilde{V}^\pi_T)\Big|\mathcal{F}_{\theta^{(n)}_{\zeta\wedge\eta}}\Big). \end{aligned} \tag{5.12}$$

By (5.4), (5.12), the Cauchy–Schwarz and Jensen inequalities,

$$I_1 \le C(f) n^{-1/4} \tag{5.13}$$

where $C(f)$ is a constant which depends only on $f$. Next, we estimate $I_2$. Recall (see Section 4 in [11]) that

$$\begin{aligned} \tilde{V}^\pi_t &= \tilde{E}^B(g(\tilde{B}_{t_1}, \dots, \tilde{B}_{t_k})|\mathcal{F}^B_t) \\ &= \tilde{V}^\pi_0 + \sum_{i=1}^{j} \int_{t_{i-1}}^{t_i} \frac{\partial q_i}{\partial x_i}(u, \tilde{B}_{t_1}, \dots, \tilde{B}_{t_{i-1}}, \tilde{B}_u)\, d\tilde{B}_u \\ &\quad + \int_{t_j}^{t} \frac{\partial q_{j+1}}{\partial x_{j+1}}(u, \tilde{B}_{t_1}, \dots, \tilde{B}_{t_j}, \tilde{B}_u)\, d\tilde{B}_u, \qquad \text{if } t \in [t_j, t_{j+1}] \quad \text{and} \\ \tilde{V}^\pi_t &= g(\tilde{B}_{t_1}, \dots, \tilde{B}_{t_k}), \qquad \text{if } t_k \le t \le T, \end{aligned} \tag{5.14}$$

where $t_0 = 0$ and the functions $q_i : [t_{i-1}, t_i] \times \mathbb{R}^i \to \mathbb{R}$ are defined inductively as follows:

$$
\begin{aligned}
&q_k(t, x_1, \dots, x_k) \\
&\quad = (2\pi(t_k - t))^{-1/2} \int_{\mathbb{R}} g(x_1, \dots, x_{k-1}, x_k + u) \exp\left(-\frac{u^2}{2(t_k - t)}\right) du \\
&\qquad \text{if } t_{k-1} \le t < t_k, q_k(t_k, x_1, \dots, x_k) = g(x_1, \dots, x_k) \quad \text{and} \\
&\text{for } i < k, \qquad q_i(t, x_1, \dots, x_i) = (2\pi(t_i - t))^{-1/2} \int_{\mathbb{R}} q_{i+1}(t_i, x_1, \dots, x_i, x_i + u) \\
&\qquad \times \exp\left(-\frac{u^2}{2(t_i - t)}\right) du \\
&\qquad \text{if } t_{i-1} \le t < t_i, q_i(t_i, x_1, \dots, x_i) = q_{i+1}(t_i, x_1, \dots, x_i, x_i).
\end{aligned}
\tag{5.15}
$$

Clearly, for any $x = (x_1, \dots, x_k)$, $y = (y_1, \dots, y_k)$ we have $|g(x) - g(y)| \le kL(f) \max_{1 \le i \le k} |x_i - y_i|$. Then it follows from (5.15) by means of the backward induction that for any $j \le k$, $|q_j(t, x_1, \dots, x_j) - q_j(t, y_1, \dots, y_j)| \le kL(f) \max_{1 \le i \le j} |x_i - y_i|$. Thus for any $j \le k$,

$$
\sup_{t \in [t_{j-1}, t_j]} \sup_{x \in \mathbb{R}^j} \left| \frac{\partial q_j}{\partial x_j}(t, x_1, \dots, x_j) \right| \le kL(f). \tag{5.16}
$$

From (5.14), (5.16) and Itô's isometry for stochastic integrals we obtain that

$$
\tilde{E}^B (\tilde{V}^\pi_{\theta^{(n)}_{\zeta \wedge \eta}} - \tilde{V}^\pi_{\sigma \wedge \theta^{(n)}_\eta})^2 \le k^2 (L(f))^2 \tilde{E}^B |\theta^{(n)}_{\zeta \wedge \eta} - \sigma \wedge \theta^{(n)}_\eta| \le k^2 (L(f))^2 \tilde{E}^B w_n,
$$

which together with (5.3) and the Cauchy–Schwarz inequality yields

$$
I_2 \le \tilde{C}(f) n^{-1/2} \tag{5.17}
$$

for some constant $\tilde{C}(f)$ which depends only on $f$. Combining (5.6), (5.10)–(5.13) and (5.17) we conclude that there is a constant $C^{(1)}(f)$ such that for any $n$, $R_n(x) - R(x) \le \epsilon + 2\delta + C^{(1)} n^{-1/4} (\ln n)^{3/4}$, and so for any initial capital $x$,

$$
R(x) \ge \limsup_{n \to \infty} R_n(x). \tag{5.18}
$$

Next we want to prove (2.26) and (2.31). Fix an initial capital $x$ and an integer $n \ge 1$. Set $(\pi, \sigma) = (\psi_n(\pi_n), \phi_n(\sigma_n))$ where $(\pi_n, \sigma_n) \in \mathcal{A}^{\xi, n}(x) \times \mathcal{T}^\xi_{0n}$ is the optimal hedge given by (3.14) and the functions $\psi_n, \phi_n$ were defined in Section 2. We can consider the portfolio $\pi = \psi_n(\pi_n)$ not only as an element in $\mathcal{A}^{B,n}(x)$ but also as an element in $\mathcal{A}^B(x)$ if we restrict the above portfolio to the interval $[0, T]$. From Lemma 3.3 we obtain that

$$
R(\pi, \sigma) - R_n(x) = R(\pi, \sigma) - R^{B,n}(\pi, \zeta_n) \tag{5.19}
$$

where, recall, $\zeta_n$ was defined in (3.15). Observe that by (2.27) and (3.14), $\sigma=\phi_n(\sigma_n)=T\wedge\theta^{(n)}_{\zeta_n}\mathbb{I}_{\zeta_n<n}+T\mathbb{I}_{\zeta_n=n}$. Since $n$ is fixed we will skip the index writing $\zeta=\zeta_n$. Given $\delta>0$ there exists a stopping time $\tau=\tau(n,\delta)$ such that

$$R(\pi,\sigma)<\delta+E^B[(Q^B(\sigma,\tau)-\tilde V^{\pi}_{\sigma\wedge\tau})^+]. \tag{5.20}$$

Let $\eta(n,\tau)=n\wedge\min\{k|\theta^{(n)}_k\geq\tau\}\in\mathcal{T}^{B,n}_{0,n}$. Observe that $\min\{k|\theta^{(n)}_k\geq\tau\}\in\mathcal{T}^{B,n}$ since $\{\min\{k|\theta^{(n)}_k\geq\tau\}\leq j\}=\{\theta^{(n)}_j\geq\tau\}\in\mathcal{F}^B_{\theta^{(n)}_j}$. From (5.9) it follows that

$$R^{B,n}(\pi,\zeta)\geq E^B[(Q^B(\theta^{(n)}_\zeta,\theta^{(n)}_\eta)-\tilde V^{\pi}_{\theta^{(n)}_{\zeta\wedge\eta}})^+]-C_1n^{-1/4}(\ln n)^{3/4}. \tag{5.21}$$

Set

$$\Gamma_1=(Q^B(\sigma,\tau)-Q^B(\theta^{(n)}_\zeta,\theta^{(n)}_\eta))^+,$$
$$\Gamma_2=|Q^B(\sigma,\tau)-Q^B(\sigma\wedge\theta^{(n)}_n,\tau\wedge\theta^{(n)}_n)|.$$

From (5.20) and (5.21) we obtain that

$$\begin{aligned}R^{B,n}(\pi,\zeta)\geq{}& E^B[(Q^B(\sigma\wedge\theta^{(n)}_n,\tau\wedge\theta^{(n)}_n)-\tilde V^{\pi}_{\theta^{(n)}_{\zeta\wedge\eta}})^+]\\ &-C_1n^{-1/4}(\ln n)^{3/4}-E^B(\Gamma_1+\Gamma_2),\\ R(\pi,\sigma)<{}&\delta+E^B[(Q^B(\sigma\wedge\theta^{(n)}_n,\tau\wedge\theta^{(n)}_n)-\tilde V^{\pi}_{\sigma\wedge\tau})^+]+E^B\Gamma_2.\end{aligned}$$

Hence,

$$\begin{aligned}R(\pi,\sigma)-R^{B,n}(\pi,\zeta)<{}&E^B[(Q^B(\sigma\wedge\theta^{(n)}_n,\tau\wedge\theta^{(n)}_n)-\tilde V^{\pi}_{\sigma\wedge\tau})^+]\\ &-E^B[(Q^B(\sigma\wedge\theta^{(n)}_n,\tau\wedge\theta^{(n)}_n)-\tilde V^{\pi}_{\theta^{(n)}_{\zeta\wedge\eta}})^+]\\ &+\delta+E^B(\Gamma_1+2\Gamma_2)+C_1n^{-1/4}(\ln n)^{3/4}.\end{aligned} \tag{5.22}$$

Observe that $\sigma\wedge\tau\wedge\theta^{(n)}_n\leq\theta^{(n)}_{\zeta\wedge\eta}$. Since $\pi\in\mathcal{A}^{B,n}(x)$, then by (2.29), $\tilde V^{\pi}_{\sigma\wedge\tau}=\tilde V^{\pi}_{\sigma\wedge\tau\wedge\theta^{(n)}_n}=\tilde E^B(\tilde V^{\pi}_{\theta^{(n)}_{\zeta\wedge\eta}}|\mathcal{F}^B_{\sigma\wedge\tau\wedge\theta^{(n)}_n})$ which together with the Jensen inequality yields that

$$\begin{aligned}&(Q^B(\sigma\wedge\theta^{(n)}_n,\tau\wedge\theta^{(n)}_n)-\tilde V^{\pi}_{\sigma\wedge\tau})^+\\ &\quad\leq\tilde E^B((Q^B(\sigma\wedge\theta^{(n)}_n,\tau\wedge\theta^{(n)}_n)-\tilde V^{\pi}_{\theta^{(n)}_{\zeta\wedge\eta}})^+|\mathcal{F}^B_{\sigma\wedge\tau\wedge\theta^{(n)}_n})\\ &\quad=E^B\Big(\frac{Z_{\sigma\wedge\tau\wedge\theta^{(n)}_n}}{Z_{\theta^{(n)}_{\zeta\wedge\eta}}}(Q^B(\sigma\wedge\theta^{(n)}_n,\tau\wedge\theta^{(n)}_n)-\tilde V^{\pi}_{\theta^{(n)}_{\zeta\wedge\eta}})^+|\mathcal{F}^B_{\sigma\wedge\tau\wedge\theta^{(n)}_n}\Big).\end{aligned} \tag{5.23}$$

Thus,

$$
\begin{aligned}
&E^B(Q^B(\sigma \wedge \theta_n^{(n)}, \tau \wedge \theta_n^{(n)}) - \tilde{V}^\pi_{\sigma\wedge\tau})^+ \\
&\quad \le E^B\biggl(\frac{Z_{\sigma\wedge\tau\wedge\theta_n^{(n)}}}{Z_{\theta^{(n)}_{\zeta\wedge\eta}}}(Q^B(\sigma \wedge \theta_n^{(n)}, \tau \wedge \theta_n^{(n)}) - \tilde{V}^\pi_{\theta^{(n)}_{\zeta\wedge\eta}})^+\biggr).
\end{aligned} \tag{5.24}
$$

By (5.22) and (5.24) we obtain that

$$
R(\pi,\sigma) - R^{B,n}(\pi,\zeta) < C_1 n^{-1/4}(\ln n)^{3/4} + \delta + E^B(\Gamma_1 + 2\Gamma_2) + I \tag{5.25}
$$

where

$$
I = E^B\biggl(\frac{Z_{\sigma\wedge\tau\wedge\theta_n^{(n)}} - Z_{\theta^{(n)}_{\zeta\wedge\eta}}}{Z_{\theta^{(n)}_{\zeta\wedge\eta}}}(Q^B(\sigma \wedge \theta_n^{(n)}, \tau \wedge \theta_n^{(n)}) - \tilde{V}^\pi_{\theta^{(n)}_{\zeta\wedge\eta}})^+\biggr).
$$

Notice that

$$
\begin{aligned}
&|\sigma\wedge\tau - \theta^{(n)}_{\zeta\wedge\eta}| \le w_n \quad \text{and} \\
&|\sigma\wedge\tau\wedge\theta_n^{(n)} - \theta^{(n)}_{\zeta\wedge\eta}| \le |\sigma\wedge\tau - \theta^{(n)}_{\zeta\wedge\eta}| \le w_n.
\end{aligned} \tag{5.26}
$$

From Itô's formula it follows that $dZ_t = \frac{\mu}{\kappa} Z_t\,dB_t + (\frac{\mu}{\kappa})^2 Z_t\,dt$, and so

$$
Z_{\theta^{(n)}_{\zeta\wedge\eta}} - Z_{\sigma\wedge\tau\wedge\theta_n^{(n)}} = \frac{\mu}{\kappa}\int_{\sigma\wedge\tau\wedge\theta_n^{(n)}}^{\theta^{(n)}_{\zeta\wedge\eta}} Z_t\,dB_t + \left(\frac{\mu}{\kappa}\right)^2 \int_{\sigma\wedge\tau\wedge\theta_n^{(n)}}^{\theta^{(n)}_{\zeta\wedge\eta}} Z_t\,dt.
$$

Set $D_n = \sup_{0\le t\le \theta_n^{(n)}\vee T} Z_t$. From (5.3), the Cauchy–Schwarz inequality and Itô's isometry we obtain that

$$
\begin{aligned}
&E^B(Z_{\theta^{(n)}_{\zeta\wedge\eta}} - Z_{\sigma\wedge\tau\wedge\theta_n^{(n)}})^2 \\
&\quad \le 2\left(\frac{\mu}{\kappa}\right)^2 E^B\int_{\sigma\wedge\tau\wedge\theta_n^{(n)}}^{\theta^{(n)}_{\zeta\wedge\eta}} Z_t^2\,dt + 2\left(\frac{\mu}{\kappa}\right)^4 E^B(w_n D_n)^2 \\
&\quad \le 2\left(\frac{\mu}{\kappa}\right)^2 E^B(w_n D_n^2) + 2\left(\frac{\mu}{\kappa}\right)^4 E^B(w_n D_n)^2 \\
&\quad \le H^{(5)} n^{-1/2}
\end{aligned} \tag{5.27}
$$

for some constant $H^{(5)}$. Since $Q^B(\sigma\wedge\theta_n^{(n)}, \tau\wedge\theta_n^{(n)}) \le A_n$ by (4.1), then by (5.27) and the Cauchy–Schwarz inequality there exists a constant $H^{(6)}$ such that

$$
I \le H^{(6)} n^{-1/4}. \tag{5.28}
$$

Next we want to estimate $E^B\Gamma_1$. Observe that if $\sigma<\tau$, then $\zeta<\eta$, and so by (5.3), (5.26) and Lemma 4.3 there exists a constant $H^{(7)}$ such that

$$\begin{aligned} E^B\Gamma_1 &\le E^B|e^{-r\sigma\wedge\tau}G_{\sigma\wedge\tau}(S^B)-e^{-r\theta^{(n)}_{\zeta\wedge\eta}}G_{\theta^{(n)}_{\zeta\wedge\beta}}(S^B)| \\ &\quad + E^B|e^{-r\sigma\wedge\tau}F_{\sigma\wedge\tau}(S^B)-e^{-r\theta^{(n)}_{\zeta\wedge\beta}}F_{\theta^{(n)}_{\zeta\wedge\eta}}(S^B)|\le H^{(7)}n^{-1/4}. \end{aligned} \tag{5.29}$$

Finally we estimate $E^B\Gamma_2$. From the definitions it follows easily that $\sigma<\tau$ is equivalent to $\sigma\wedge\theta^{(n)}_n<\tau\wedge\theta^{(n)}_n$, and so from (5.26) it follows that there exists a constant $H^{(8)}$ such that

$$\begin{aligned} E^B\Gamma_2 &\le E^B|e^{-r\sigma\wedge\tau}G_{\sigma\wedge\tau}(S^B)-e^{-r\theta^{(n)}_n\wedge\sigma\wedge\tau}G_{\theta^{(n)}_n\wedge\sigma\wedge\tau}(S^B)| \\ &\quad + E^B|e^{-r\sigma\wedge\tilde\tau}F_{\sigma\wedge\tilde\tau}(S^B)-e^{-r\theta^{(n)}_n\wedge\sigma\wedge\tau}F_{\theta^{(n)}_n\wedge\sigma\wedge\tau}(S^B)| \\ &\le H^{(8)}n^{-1/4}. \end{aligned} \tag{5.30}$$

Since $\delta$ is arbitrary, then combining (5.19), (5.25) and (5.28)–(5.30) we conclude that there is a constant $C^{(2)}$ (which does not depend on the initial capital $x$) such that $R(\pi,\sigma)-R_n(x)\le C^{(2)}n^{-1/4}(\ln n)^{3/4}$. Recall that $(\pi,\sigma)=(\psi_n(\pi_n),\phi_n(\sigma_n))$, and so for all $n\ge 1$,

$$R(\psi_n(\pi_n),\phi_n(\sigma_n))-R_n(x)\le C^{(2)}n^{-1/4}(\ln n)^{3/4}, \tag{5.31}$$

which together with (5.18) completes the proof of Theorems 2.1 and 2.2.

**6. Additional estimates for American options.** In the case of American options in BS markets the definitions (2.12) of the shortfall risks take on the following form:

$$\begin{aligned} R(\pi) &= \sup_{\tau\in\mathcal{T}^B_{0T}} E^B[(\tilde Y_\tau-\tilde V^\pi_\tau)^+], \qquad \pi\in\mathcal{A}^B \quad \text{and} \\ R(x) &= \inf_{\pi\in\mathcal{A}^B(x)} R(\pi) \end{aligned} \tag{6.1}$$

where $\tilde Y_t$ is defined after (2.10). Similarly for CRR models we have

$$\begin{aligned} R_n(\pi) &= \max_{\tau\in\mathcal{T}^\xi_{0n}} E^\xi_n[(\tilde Y^{(n)}_\tau-\tilde V^\pi_\tau)^+], \qquad \pi\in\mathcal{A}^{\xi,n} \quad \text{and} \\ R_n(x) &= \inf_{\pi\in\mathcal{A}^{\xi,n}(x)} R_n(\pi). \end{aligned} \tag{6.2}$$

THEOREM 6.1. *There exists a constant $\mathcal{C}$ such that for any initial capital $x$ and $n\in\mathbb{N}$ in addition to (2.26) we have*

$$R_n(x)\le R(x)+\mathcal{C}n^{-1/4}(\ln n)^{3/4}. \tag{6.3}$$

REMARK 6.2. It is easy to see that all proofs of previous sections go through for American options simplifying the corresponding arguments. Namely, assume formally in previous sections that the seller is allowed to stop only at time $T$ in the continuous-time case and at time $n$ at the $n$-step CRR model; then since $\phi_n(n) = T$ [see (2.27)] all proofs above will go through and we derive the results of Section 2 for corresponding American options, as well, assuming (2.1)–(2.2) for payoffs. In general, American options can be considered as partial cases of game options where penalties are chosen so high that it will not be wise for the seller to stop until the expiration time; but in order to apply our results from previous sections to such game options directly we have to construct such penalties satisfying conditions (2.1)–(2.2), which is not very easy.

The dynamical programming algorithm that we used in order to calculate optimal hedges for Israeli options is also valid in the American options case. Namely, similarly to (3.12)–(3.13) define

$$\begin{aligned} &J_n^n(y, u_1, u_2, \ldots, u_n) \\ &\quad = (f_n^n(u_1, \ldots, u_n) - y)^+, \\ &J_k^n(y, u_1, \ldots, u_k) \\ &\quad = \max\Bigg((f_k^n(u_1, \ldots, u_k) - y)^+, \\ &\qquad \inf_{u \in I_n(y)} \Bigg[p^{(n)} J_{k+1}^n\Bigg(y + ua_1^{(n)}, u_1, \ldots, u_k, \sqrt{\frac{T}{n}}\Bigg) \\ &\qquad\qquad + (1 - p^{(n)}) J_{k+1}^n\Bigg(y + ua_2^{(n)}, u_1, \ldots, u_k, -\sqrt{\frac{T}{n}}\Bigg)\Bigg]\Bigg) \\ &\qquad\qquad\qquad \text{for } k = n-1, n-2, \ldots, 0 \end{aligned} \tag{6.4}$$

and

$$\begin{aligned} &h_k^n(y, x_1, \ldots, x_k) \\ &\quad = \operatorname*{arg\,min}_{u \in I_n(y)} \Bigg[p^{(n)} J_{k+1}^n\Bigg(y + ua_1^{(n)}, u_1, \ldots, u_k, \sqrt{\frac{T}{n}}\Bigg) \\ &\qquad\qquad + (1 - p^{(n)}) J_{k+1}^n\Bigg(y + ua_2^{(n)}, u_1, \ldots, u_k, -\sqrt{\frac{T}{n}}\Bigg)\Bigg], \\ &\qquad\qquad\qquad k < n. \end{aligned} \tag{6.5}$$

Similarly to (3.14), for a given initial capital $x$ and $n \in \mathbb{N}$ define an *admissible* self- financing strategy $\pi_n$ by

$$\tilde{V}_0^{\pi_n} = x \quad \text{and}$$

$$\tilde V^{\pi_n}_{k+1}=\tilde V^{\pi_n}_k+h^n_k(\tilde V^{\pi_n}_k,e^{\kappa\sqrt{T/n}\xi_1},\dots,e^{\kappa\sqrt{T/n}\xi_k})(e^{\kappa\sqrt{T/n}\xi_{k+1}}-1)\quad\text{for } k>0. \tag{6.6}$$

As in Lemma 3.3 we have that

$$R_n(\pi_n)=R_n(x). \tag{6.7}$$

For American options we can also improve Theorem 2.2 as follows.

THEOREM 6.3. *For any $n$ let $\pi_n\in\mathcal{A}^{\xi,n}(x)$ be the optimal hedge constructed in (6.6); then*

$$\lim_{n\to\infty}R(\psi_n(\pi_n))=R(x). \tag{6.8}$$

*Furthermore, there exists a constant $\tilde{\mathcal{C}}$ such that*

$$R(\psi_n(\pi_n))\le R(x)+\tilde{C}n^{-1/4}(\ln n)^{3/4}. \tag{6.9}$$

In order to derive these results we will need several lemmas. Let $n\in\mathbb{N}$ and consider the restriction of the measures $P^B,\tilde P^B$ to the $\sigma$-algebra $\mathcal{G}^{B,n}_n$. Set $W_n=\frac{dP^B}{d\tilde P^B}|\mathcal{G}^{B,n}_n$. Observe that $\int_A W_n\,d\tilde P^B=P^B(A)$ for any $A\in\mathcal{G}^{B,n}_n$. Since $A\in\mathcal{F}^B_{\theta^{(n)}_n}$, then $\int_A Z_{\theta^{(n)}_n}\,d\tilde P^B=P^B(A)$, and so

$$W_n=\tilde E^B(Z_{\theta^{(n)}_n}|\mathcal{G}^{B,n}_n). \tag{6.10}$$

LEMMA 6.4. *There exists a constant $C_2$ such that for any $n$,*

$$\tilde E^B(W_n-Z_{\theta^{(n)}_n})^2\le C_2n^{-1/2}. \tag{6.11}$$

PROOF. We know that $Z_{\theta^{(n)}_n}=\exp(aB^*_{\theta^{(n)}_n}+b\theta^{(n)}_n)$ where $a=\frac{\mu}{\kappa}$ and $b=-\frac{\mu}{2}-\frac{\mu^2}{2\kappa^2}$. Set $V_n=\exp(aB^*_{\theta^{(n)}_n}+bT)$ which is clearly $\mathcal{G}^{B,n}_n$-measurable. Since conditional expectation is an orthogonal projection it follows from (6.10) that

$$\tilde E^B(W_n-Z_{\theta^{(n)}_n})^2\le\tilde E^B(V_n-Z_{\theta^{(n)}_n})^2. \tag{6.12}$$

Using Cauchy–Schwarz and Chebyshev inequalities together with the inequality $|e^{bx}-1|\le|b|e^{|b|}|x|$ for $-1\le x\le 1$ we obtain

$$\begin{aligned}&\tilde E^B(V_n-Z_{\theta^{(n)}_n})^2\\&\quad\le\tilde E^B[\mathbb{I}_{\{1<|\theta^{(n)}_n-T|\}}(V_n^2+Z^2_{\theta^{(n)}_n})]\end{aligned}$$

$$
\begin{aligned}
&+\tilde{E}^B[\mathbb{I}_{\{1\geq|\theta_n^{(n)}-T|\}}V_n^2|e^{|b|(\theta_n^{(n)}-T)}-1|^2]\\
&\leq(\tilde{E}^B(V_n^2+Z_{\theta_n^{(n)}}^2)^2)^{1/2}(\tilde{E}^B\mathbb{I}_{\{1<|\theta_n^{(n)}-T|\}})^{1/2}\\
&\quad+\tilde{E}^B(b^2e^{2|b|}V_n^2|\theta_n^{(n)}-T|^2)\\
&\leq(\tilde{E}^B(V_n^2+Z_{\theta_n^{(n)}}^2)^2)^{1/2}(\tilde{E}^B|\theta_n^{(n)}-T|^2)^{1/2}\\
&\quad+b^2e^{2|b|}(\tilde{E}^BV_n^4)^{1/2}(\tilde{E}^B|\theta_n^{(n)}-T|^4)^{1/2}\\
&\leq C_2n^{-1/2}
\end{aligned}\tag{6.13}
$$

for some constant $C_2$. Now (6.11) follows from (6.12) and (6.13), completing the proof. □

Lemma 6.5. *For $n\in\mathbb{N}$ let $\{M_i\}_{i=0}^n$ be a martingale with respect to the filtration $\{\mathcal{F}^B_{\theta_i^{(n)}}\}_{i=0}^n$ and the measure $\tilde{P}^B$. Set $\tilde{M}_i=\tilde{E}^B(M_i|\mathcal{G}_n^{B,n})$. Then $\{\tilde{M}_i\}_{i=0}^n$ is a martingale with respect to the filtration $\{\mathcal{G}_i^{B,n}\}_{i=0}^n$ and the measure $\tilde{P}^B$.*

Proof. For a fixed $0\leq k\leq n$ set $\Psi=M_k$, $\mathcal{F}=\mathcal{G}_k^{B,n}$, $\mathcal{K}=\sigma(B^*_{\theta_{k+1}^{(n)}}-B^*_{\theta_k^{(n)}},\dots,B^*_{\theta_n^{(n)}}-B^*_{\theta_{n-1}^{(n)}})$ and $\mathcal{H}=\mathcal{G}_n^{B,n}$. Using Remark 4.3 in [7] we obtain

$$
\begin{aligned}
\tilde{M}_k&=\tilde{E}^B(M_k|\mathcal{G}_n^{B,n})=\tilde{E}^B(M_k|\mathcal{G}_k^{B,n})=\tilde{E}^B(\tilde{E}^B(M_n|\mathcal{F}^B_{\theta_k^{(n)}})|\mathcal{G}_k^{B,n})\\
&=\tilde{E}^B(M_n|\mathcal{G}_k^{B,n})=\tilde{E}^B(\tilde{E}^B(M_n|\mathcal{G}_n^{B,n})|\mathcal{G}_k^{B,n})=\tilde{E}^B(\tilde{M}_n|\mathcal{G}_k^{B,n})
\end{aligned}\tag{6.14}
$$

and the result follows. □

Next, we will need some additional estimates. For any initial capital $x$ and $n\in\mathbb{N}$ define

$$
J_n(x)=\inf_{\pi\in\mathcal{A}^B(x)}\sup_{\tau\in\mathcal{T}_{0,n}^{B,n}}E^B[(\tilde{Y}^{B,n}_{\tau T/n}-\tilde{V}^\pi_{T\wedge\theta_\tau^{(n)}})^+]\tag{6.15}
$$

where, recall, $\tilde{Y}_t^{B,n}$ is defined after (3.6). The following inequality is the main point which we cannot extend directly to game options in view of the additional infimum in stopping times of the option seller there.

Lemma 6.6. *There exists a constant $C_3$ such that for any initial capital $x$ and $n\in\mathbb{N}$,*

$$
J_n(x)\leq R(x)+C_3n^{-1/4}.\tag{6.16}
$$

PROOF. Fix $n\in\mathbb{N}$ and an initial capital $x$. By using (5.9) for $\eta=n$ we get that $\sup_{\zeta\in\mathcal{T}_{0,n}^{B,n}} E^B|\tilde{Y}_{\theta_\zeta^{(n)}}-\tilde{Y}^{B,n}_{\zeta T/n}|\le C_1 n^{-1/4}(\ln n)^{3/4}$. From (5.3) and Lemma 4.3 it follows that there exists a constant $\tilde{C}_1$ such that $\sup_{\zeta\in\mathcal{T}_{0,n}^{B,n}} E^B|\tilde{Y}_{\theta_\zeta^{(n)}}-\tilde{Y}_{T\wedge\theta_\zeta^{(n)}}|\le\tilde{C}_1 n^{-1/4}$. Thus for $C_3=C_1+\tilde{C}_1$, $\sup_{\zeta\in\mathcal{T}_{0,n}^{B,n}} E^B|\tilde{Y}_{T\wedge\theta_\zeta^{(n)}}-\tilde{Y}^{B,n}_{\zeta T/n}|\le C_3 n^{-1/4}$. Hence,

$$\begin{aligned}J_n(x)&=\inf_{\pi\in\mathcal{A}^B(x)}\sup_{\zeta\in\mathcal{T}_{0,n}^{B,n}} E^B[(\tilde{Y}^{B,n}_{\zeta T/n}-\tilde{V}^\pi_{T\wedge\theta_\zeta^{(n)}})^+]\\&\le C_3n^{-1/4}+\inf_{\pi\in\mathcal{A}^B(x)}\sup_{\zeta\in\mathcal{T}_{0,n}^{B,n}} E^B[(\tilde{Y}_{T\wedge\theta_\zeta^{(n)}}-\tilde{V}^\pi_{T\wedge\theta_\zeta^{(n)}})^+]\\&\le C_3n^{-1/4}+\inf_{\pi\in\mathcal{A}^B(x)}\sup_{\tau\in\mathcal{T}_{0,T}^{B}} E^B[(\tilde{Y}_\tau-\tilde{V}^\pi_\tau)^+]\\&=C_3n^{-1/4}+R(x).\end{aligned}$$ □

For any initial capital $x$ and $n\in\mathbb{N}$ define

$$E_n(x)=\inf_{\pi\in\mathcal{A}^B(x)}\sup_{\tau\in\mathcal{T}_{0,n}^{B,n}}\tilde{E}^B[(\tilde{Y}^{B,n}_{\tau T/n}-\tilde{V}^\pi_{T\wedge\theta_\tau^{(n)}})^+W_n]\tag{6.17}$$

where, recall, $W_n$ is defined in (6.10). From (6.15) it follows that $J_n(x)=\inf_{\pi\in\mathcal{A}^B(x)}\sup_{\tau\in\mathcal{T}_{0,n}^{B,n}}\tilde{E}^B[(\tilde{Y}^{B,n}_{\tau T/n}-\tilde{V}^\pi_{T\wedge\theta_\tau^{(n)}})^+Z_{\theta_n^{(n)}}]$. Thus from (6.11) and the Cauchy–Schwarz inequality we obtain

$$\begin{aligned}&|E_n(x)-J_n(x)|\\&\quad\le\sup_{\zeta\in\mathcal{T}_{0,n}^{B,n}}\tilde{E}^B[|W_n-Z_{\theta_n^{(n)}}|\tilde{Y}^{B,n}_{\zeta T/n}]\\&\quad\le(\tilde{E}^B(W_n-Z_{\theta_n^{(n)}})^2)^{1/2}\sup_{\zeta\in\mathcal{T}_{0,n}^{B,n}}(\tilde{E}^B(\tilde{Y}^{B,n}_{\zeta T/n})^2)^{1/2}\\&\quad\le C_4n^{-1/4},\end{aligned}$$

for some constant $C_4$. This together with Lemma 6.3 yields that there exists a constant $C_5$ such that

$$E_n(x)\le R(x)+C_5n^{-1/4}(\ln n)^{3/4}.\tag{6.18}$$

Now we return to the proof of Theorems 6.1 and 6.3. Fix an initial capital $x$ and $n\in\mathbb{N}$. Analogously to (3.7) define

$$R^{B,n}(x)=\inf_{\pi\in\mathcal{A}^{B,n}(x)}\sup_{\tau\in\mathcal{T}_{0,n}^{B,n}}E^B[(\tilde{Y}^{B,n}_{\tau T/n}-\tilde{V}^\pi_{\theta_\tau^{(n)}})^+],\tag{6.19}$$

where $\mathcal{A}^{B,n}(x)$ is defined in (2.29). Similarly to Lemma 3.3

$$R_n(x) = R^{B,n}(x). \tag{6.20}$$

Choose $\epsilon > 0$. There exists $\pi \in \mathcal{A}^{B,M}(x)$ such that

$$\sup_{\tau \in \mathcal{T}_{0,n}^{B,n}} \tilde{E}^B[(\tilde{Y}^{B,n}_{\tau T/n} - \tilde{V}^{\pi}_{T \wedge \theta^{(n)}_{\tau}})^+ W_n] < E_n(x) + \epsilon. \tag{6.21}$$

The sequence $\{\tilde{V}^{\pi}_{T\wedge\theta^{(n)}_k}\}_{k=0}^n$ is a martingale with respect to the filtration $\{\mathcal{F}^B_{\theta^{(n)}_k}\}_{k=0}^n$ and the martingale measure $\tilde{P}^B$. Define

$$\tilde{M}_k = \tilde{E}^B(\tilde{V}^{\pi}_{T\wedge\theta^{(n)}_k} | \mathcal{G}^{B,n}_n), \qquad 0 \le k \le n. \tag{6.22}$$

From Lemma 6.5 it follows that $\{\tilde{M}_k\}_{k=0}^n$ is a martingale with respect to the filtration $\{\mathcal{G}^{B,n}_k\}_{k=0}^n$ and the measure $\tilde{P}^B$. Thus for any $k \le n$ there exists a measurable function $f_k : \{-\sqrt{\frac{T}{n}}, \sqrt{\frac{T}{n}}\}^k \to \mathbb{R}_+$ such that $\tilde{M}_k = f_k(B^*_{\theta^{(n)}_1}, \dots, B^*_{\theta^{(n)}_k} - B^*_{\theta^{(n)}_{k-1}})$. Thus the sequence $\{f_k(\sqrt{\frac{T}{n}}\xi_1, \dots, \sqrt{\frac{T}{n}}\xi_k)\}_{k=0}^n$ is a martingale with respect to the filtration $\{\mathcal{F}^\xi_k\}_{k=0}^n$ and the measure $\tilde{P}^\xi_n$. Since the CRR markets are complete and $\tilde{M}_0 = x$, it follows that there exists a portfolio $\pi^\xi \in \mathcal{A}^{\xi,n}(x)$ such that for any $k \le n$ $\tilde{V}^{\pi^\xi}_k = f_k(\sqrt{\frac{T}{n}}\xi_1, \dots, \sqrt{\frac{T}{n}}\xi_k)$. Hence, we obtain for the portfolio $\tilde{\pi} = \psi_n(\pi^\xi) \in \mathcal{A}^{B,n}(x)$ that for any $k \le n$,

$$\tilde{V}^{\tilde{\pi}}_{\theta^{(n)}_k} = \tilde{M}_k. \tag{6.23}$$

Thus by (6.19)–(6.20),

$$\begin{aligned} R_n(x) &\le \sup_{\zeta \in \mathcal{T}^{B,n}_{0,n}} E^B[(\tilde{Y}^{B,n}_{\zeta T/n} - \tilde{V}^{\tilde{\pi}}_{\theta^{(n)}_\zeta})^+] \\ &= \sup_{\zeta \in \mathcal{S}^{B,n}_{0,n}} E^B[(\tilde{Y}^{B,n}_{\zeta T/n} - \tilde{V}^{\tilde{\pi}}_{\theta^{(n)}_\zeta})^+] \end{aligned} \tag{6.24}$$

where the last equality follows from the fact that $(\tilde{Y}^{B,n}_{kT/n} - \tilde{V}^{\tilde{\pi}}_{\theta^{(n)}_k})^+$ is $\mathcal{G}^{B,n}_k$-measurable (for any $k$). Since $W_n$ is $\mathcal{G}^{B,n}_n$-measurable, then from (6.22) and (6.23) it follows that for any $\zeta \in \mathcal{S}^{B,n}_{0,n}$,

$$W_n(\tilde{Y}^{B,n}_{\zeta T/n} - \tilde{V}^{\tilde{\pi}}_{\theta^{(n)}_\zeta}) = \tilde{E}^B(W_n(\tilde{Y}^{B,n}_{\zeta T/n} - \tilde{V}^{\pi}_{T\wedge\theta^{(n)}_\zeta}) | \mathcal{G}^{B,n}_n), \tag{6.25}$$

and from Jensen's inequality we obtain that

$$\tilde{E}^B[W_n(\tilde{Y}^{B,n}_{\zeta T/n} - \tilde{V}^{\tilde{\pi}}_{\theta^{(n)}_\zeta})^+] \le \tilde{E}^B[W_n(\tilde{Y}^{B,n}_{\zeta T/n} - \tilde{V}^{\pi}_{T\wedge\theta^{(n)}_\zeta})^+]. \tag{6.26}$$

By (6.21), (6.24), (6.26) and the definition of $W_n$,

$$\begin{aligned} R_n(x) &\le \sup_{\zeta\in\mathcal{S}^{B,n}_{0,n}} \tilde{E}^B[(\tilde{Y}^{B,n}_{\zeta T/n} - \tilde{V}^{\tilde{\pi}}_{\theta^{(n)}_\zeta})^+ W_n] \\ &\le \sup_{\zeta\in\mathcal{S}^{B,n}_{0,n}} \tilde{E}^B[(\tilde{Y}^{B,n}_{\zeta T/n} - \tilde{V}^{\pi}_{T\wedge\theta^{(n)}_\zeta})^+ W_n] \qquad (6.27) \\ &< E_n(x) + \epsilon. \end{aligned}$$

Since $\epsilon > 0$ is arbitrary, then $R_n(x) \le E_n(x)$ which together with (6.18) completes the proof of Theorem 6.1. Using the inequality (5.31) for the case of American options it follows that for any $n$,

$$R(\psi_n(\pi_n)) - R_n(x) \le C^{(2)} n^{-1/4} (\ln n)^{3/4}, \qquad (6.28)$$

which together with Theorem 6.1 completes the proof of Theorem 6.3.

Institute of Mathematics
Hebrew University
Jerusalem 91904
Israel
E-mail: yann1@math.huji.ac.il
kifer@math.huji.ac.il